\documentclass[11pt,a4paper,leqno]{amsart}

\usepackage{amsfonts,amsmath,amssymb,amsthm}
\usepackage{hyperref}

\def\C{{\mathbb C}}
\def\R{{\mathbb R}}
\def\N{{\mathbb N}}

\def\O{{\mathcal O}}

\def\virgp{\raise 2pt\hbox{,}}

\def\ul{{\underline u}}
\def\eps{{\varepsilon}}
\def\g{\gamma}
\def\om{\omega}

\def\sumetage#1#2{
\sum_{\scriptstyle {#1}\atop\scriptstyle {#2}}
}

\def\({\left(}
\def\){\right)}
\def\<{\left\langle}
\def\>{\right\rangle}
\def\le{\leqslant}
\def\ge{\geqslant}

\def\d{{\partial}}

\DeclareMathOperator{\RE}{Re}
\DeclareMathOperator{\IM}{Im}

\theoremstyle{plain}
\newtheorem{theorem}{Theorem}[section]
\newtheorem{definition}[theorem]{Definition}
\newtheorem{lemma}[theorem]{Lemma}

\newtheorem{proposition}[theorem]{Proposition}
\theoremstyle{remark}
\newtheorem*{H1}{Assumption ($H1$)}
\newtheorem*{H2}{Assumption ($H2$)}
\newtheorem*{H2tilde}{Assumption ($\widetilde{H2}$)}
\newtheorem{remark}[theorem]{Remark}

\newtheorem{example}[theorem]{Example}

\numberwithin{equation}{section}

\begin{document}

\title[Analyticity of the scattering operator]{Analyticity of the
  scattering operator for 
semilinear dispersive equations} 
\author[R. Carles]{R{\'e}mi Carles}
\address[R. Carles]{CNRS \& Universit{\'e} Montpellier 2 \\
Math{\'e}matiques, CC 051\\
Place Eug{\`e}ne Bataillon\\
34095 Montpellier cedex 5\\
France}
\email{Remi.Carles@math.cnrs.fr}
\author[I. Gallagher]{Isabelle Gallagher}
\address[I. Gallagher]{Institut de Math{\'e}matiques de Jussieu UMR 7586\\
 Universit{\'e} Paris VII\\
175, rue du Chevaleret\\ 
75013 Paris\\
France}
\email{Isabelle.Gallagher@math.jussieu.fr}
\thanks{This work was partially supported by the ANR project SCASEN}
\begin{abstract}
We present a general algorithm to show that a scattering operator
associated to a semilinear dispersive equation is real
analytic, and to compute the coefficients of its Taylor series at any
point. We illustrate this method in the case of the Schr\"odinger
equation with power-like nonlinearity or with Hartree type
nonlinearity, and in the case of the wave and Klein--Gordon equations
with power nonlinearity. Finally, we discuss the link of this approach
with inverse scattering, and with complete integrability.  
\end{abstract}
\subjclass[2000]{Primary: 35P25. Secondary: 35B30; 35C20; 35L05; 35Q55; 37K10}
\keywords{Schr\"odinger equation; wave equation; Klein--Gordon equation; scattering operator; inverse scattering; integrability}
\maketitle

\section{Introduction}
\label{sec:intro}
The local and global well-posedness of semilinear dispersive equations
has attracted a lot of attention for the past years. In general, when
global well-posedness is established, the existence of a scattering
operator, comparing the nonlinear dynamics and the linear one, is a
rather direct by-product. Unlike in the linear case (see
e.g. \cite{LaxPhillips,ReedSimon3,Yafaev}), besides continuity, very
few properties of these 
nonlinear scattering operators are known. A first natural question,
which can be found in \cite[pp.~121--122]{Reed76}, consists in
investigating the real analyticity of the scattering operators. A
positive answer is available in some very specific cases: see
\cite{BaezZhou89,BaezZhou90,Kumlin} for the cubic wave and
Klein--Gordon equation in 3D, and \cite{R3Analytic} for the Hartree
equation in 3D. In this paper, we extend these results to a more
general class of dispersive equations, including the nonlinear
Schr\"odinger equation and the nonlinear wave equation, in space
dimension $n\le 4$ (such an assumption is needed for the power
nonlinearity to be both analytic and energy-subcritical or
critical). Moreover, unlike in 
\cite{BaezZhou89,BaezZhou90,Kumlin,R3Analytic}, we do not use an
abstract analytic implicit function theorem: we construct directly the
terms of the series \emph{via} a general abstract lemma, thus
extending the approach of S.~Masaki 
\cite{Masaki07}. We then show that the series is converging, working
in suitable spaces based on dispersive properties provided by
Strichartz estimates. In general, these estimates are a direct
by-product of the proof of the existence of a nonlinear scattering
operator. 
\smallbreak

Before being more precise about the results presented here, we briefly
recall the approach for (short range) scattering theory in the context
of semilinear 
dispersive equations. The main examples we have in mind are the
nonlinear Schr\"odinger equation
\begin{equation}
  \label{eq:NLSintro}
  i\d_t u+\frac{1}{2}\Delta u = \lambda \lvert u\rvert^{p-1}u,\quad (t,x)\in
  \R\times \R^n ,
\end{equation}
the Hartree equation
\begin{equation}
  \label{eq:r3intro}
  i\d_t u+\frac{1}{2}\Delta u = \lambda\(\lvert x\rvert^{-\gamma}\ast \lvert
  u\rvert^2\)u,\quad (t,x)\in 
  \R\times \R^n ,
\end{equation}
the nonlinear wave and Klein--Gordon  equations
\begin{eqnarray} 
  \d_t^2 u -\Delta u  + \lambda u^p=0, \quad (t,x)\in
  \R\times \R^n \label{eq:waveintro}\\
 \d_t^2 u -\Delta u +  u + \lambda u^p=0, \quad (t,x)\in
  \R\times \R^n .\label{eq:nklgintro}
\end{eqnarray}
Up to considering the unknown $(u,\d_t u)$ instead of $u$ alone in
\eqref{eq:waveintro}, \eqref{eq:nklgintro}, Duhamel's formula reads, in
all these examples, 
\begin{equation}
  \label{eq:Duhamelintro}
  u(t) = U(t)u_0 +\int_{t_0}^t U(t-s)\(F(u(s)\)ds,
\end{equation}
where $U(\cdot)$ is the group associated to the linear equation
($\lambda=0$), and $t_0$ corresponds to the time for which initial
data are prescribed:
\begin{equation}\label{eq:CIintro}
  U(-t)u(t)\big|_{t=t_0} = u_0.
\end{equation}
In the study of the Cauchy problem, one usually considers the case
$t_0=0$. In scattering theory, the first standard step consists in
solving the Cauchy problem near infinite time: $t_0=\pm \infty$. To
consider forward in time propagation, assume $t_0=-\infty$. To define
the wave operator $W_-$, one has to solve the Cauchy problem
\eqref{eq:Duhamelintro}-\eqref{eq:CIintro} with $t_0=-\infty$, on some
time interval of the form $]-\infty,T]$, for some finite
$T$. Classically, this step is achieved by a fixed point argument in
suitable function spaces. This may yield a time $T\ll -1$, that is,
``close'' to $-\infty$ (but finite). Suppose that the classical Cauchy
problem enables us to define $u$ up to time $t=0$. Then the wave
operator $W_-$ is defined by
\begin{equation*}
  W_- u_0= u_{\mid t=0} .
\end{equation*}
The second step consists in inverting the wave operators. For initial
data prescribed at time $t=0$, suppose that we can construct a
solution which is defined globally in time (or in the future only, for
our purpose). Inverting the wave operators (that is, proving the
asymptotic completeness) consists in showing that
nonlinear effects become negligible for large time, and that we can
find $u_+$ such that $u(t)\sim U(t)u_+$ as $t\to +\infty$:
\begin{equation*}
  u_+ = W_+^{-1}u_{\mid t=0}. 
\end{equation*}
The scattering operator $S$ is then defined by
\begin{equation*}
  S u_0 = W_+^{-1} W_- u_0 = u_+. 
\end{equation*}
In general, for small data, the scattering operator $S$ can be
constructed in one step only, thanks to a bootstrap argument in spaces
based on Strichartz estimates. For large data, one must expect $T\ll
-1$ in general. The solution is then made global thanks to \emph{a
  priori} estimates, such as the conservation of a positive energy
($\lambda>0$ in the above examples). The proof of asymptotic
completeness usually relies on different arguments: Morawetz
estimates, or existence of an extra evolution law (\emph{e.g.}
pseudo-conformal evolution law). In many cases, these arguments make
it possible to define the scattering operator. The continuity of this
operator is usually an easy consequence of its construction (provided
that the proof does not rely on compactness arguments). Finer
properties, such as real analyticity, are not straightforward. We
emphasize again that contrary to the 
case of the wave operators,
real analyticity of the scattering operator (for arbitrary data)
cannot be a mere consequence of the fixed point method used to
construct solutions; we show here 
that real analyticity of the scattering operator is 
very often a  consequence of the (global in time) estimates which are
established in order to show that there is scattering. In all this
paper, by ``analytic'', we mean ``real analytic'':
\begin{definition}
  Let $X$ and $Y$ be Banach spaces, and consider an operator $A:X\to
  Y$. We say that $A$ is real analytic (or simply analytic) from $X$
  to~$Y$ if~$A$ is
  infinitely Fr\'echet-differentiable at every point of $X$, with a
  locally norm-convergent series: for all $f\in X$, there exists
  $\eps_0>0$, such that for all $g\in X$, $\|g\|_X \le  1$, we can
  find $(w_j)_{j\in 
  \N}\in Y^{\N}$ such that for $0<\eps\le \eps_0$,
\begin{equation*}
  \sum_{j=0}^\infty \eps^j
  \lVert w_j\rVert_Y<\infty,\quad \text{and}\quad A(f+\eps g) = A(f) + \eps
  \displaystyle 
  \sum_{j=0}^\infty \eps^j w_j.
\end{equation*}
\end{definition}
First, it should be noted that the analyticity of scattering operators
near the origin can be obtained rather directly in general, by
applying a fixed point argument with analytic parameters. Of
course, if the nonlinearity is not analytic, one must not expect the
scattering operator to be analytic. As an illustration, consider the
nonlinear Schr\"odinger equation \eqref{eq:NLSintro}. 
As noticed in \cite{CaWigner} (in the case $n=1$), and following the
approach of \cite{PG96}, the first terms of the asymptotic expansion
of the nonlinear scattering operator $S$ near the origin are given by:
\begin{equation*}
  S\(\eps u_-\) = \eps u_- -i\eps^{p}\int_{-\infty}^{+\infty} 
e^{-i\frac{t}{2}\Delta}\(\left\lvert
  e^{i\frac{t}{2}\Delta}u_-\right\rvert^{p-1}e^{i\frac{t}{2}\Delta}u_-\)
dt + \O_{L^2}\(\eps^{2p-1}\).
\end{equation*}
The complete proof of this relation is available in \cite{COMRL} in the
$L^2$-critical case $p=1+4/n$, for any $n\ge 1$. This shows that if
$p$ is not an integer, the operator $S$ is not analytic near the
origin: it is H\"older continuous, of order~$p$ and not better. We
shall therefore consider only analytic nonlinearities:
in~\eqref{eq:waveintro}, \eqref{eq:nklgintro}, we shall always assume 
that $p$ is an integer, and in~\eqref{eq:NLSintro}, we shall assume
that $p$ is an odd integer. We can now state two typical results of
our approach. Denote
\begin{align*}
  \Sigma &= \{ f\in H^1(\R^n),\quad x\mapsto \lvert x\rvert f(x)\in
  L^2(\R^n)\}
.
\end{align*}
This space is naturally a Hilbert space. The main results of the
paper are the following. 
\begin{theorem}\label{theo:NLS}
  Let $1\le n\le 4$ and $\lambda>0$. Assume that $p\ge 3$ is an odd
  integer, with in addition
  \begin{itemize}
  \item $p\ge 5$ if $n=1$.
\item $p=3$ or $5$ if $n=3$.
\item $p=3$ if $n=4$. 
  \end{itemize}
Then the wave and scattering operators associated to the nonlinear
Schr\"o\-ding\-er equation 
\eqref{eq:NLSintro} are analytic from~$\Sigma$ to $\Sigma$. If
moreover $p\ge 7$ for $n=1$ and $p\ge 5$ for $n=2$, then the wave and
scattering operators associated to 
\eqref{eq:NLSintro} are analytic from~$H^1(\R^n)$ to $H^1(\R^n)$.
\end{theorem}
\begin{theorem}\label{theo:r3}
  Let $n\ge 3$ and $\lambda>0$. Assume that $2\le \g<\min(4,n)$. 
Then the wave and scattering operators associated to the Hartree equation
\eqref{eq:r3intro} are analytic from~$\Sigma$ to $\Sigma$. If
moreover $\g>2$, then the wave and scattering operators
associated to 
\eqref{eq:r3intro} are analytic from~$H^1(\R^n)$ to $H^1(\R^n)$.
\end{theorem}
\begin{theorem}\label{theo:NLW}
  Let $\lambda>0$. Assume that either $(n,p)=(3,5)$ or $(n,p)=(4,3)$. 
Then the wave and scattering operators associated 
to the nonlinear wave equation~\eqref{eq:waveintro} are analytic
$\dot H^1(\R^n)\times L^2(\R^n)$ to
$\dot H^1(\R^n)\times L^2(\R^n)$.
\end{theorem}
\begin{theorem}\label{theo:NLKG}
  Let $1\le n\le 4$ and $\lambda>0$. Assume that $p\ge 3$ is an odd
  integer, with   
\begin{itemize}
  \item $p\ge 7$ if $n=1$.
\item $p\ge 5$ if $n=2$. 
\item $p=3$ or $5$ if $n=3$.
\item $p=3$ if $n=4$. 
  \end{itemize}
The wave and scattering operators associated 
to the nonlinear Klein--Gordon equation~\eqref{eq:nklgintro} are
analytic from~$H^1(\R^n)\times L^2(\R^n)$ to
$H^1(\R^n)\times L^2(\R^n)$.
\end{theorem}
 


\subsubsection*{Notation} If~$A$ and~$B$ are two real numbers, we will
write~$A  \lesssim B$ if there is 
a universal constant~$C$, which does not depend on varying parameters
of the problem, such that~$A \le CB$. If~$A  \lesssim B$ 
and~$B  \lesssim A$, then we will write~$A \sim B$.

\section{An abstract result}
\label{sec:abstract}
In this section we intend to study an abstract semilinear equation,
and to present the assumptions we will make in order to conclude to  
the analyticity of the nonlinear scattering operator associated to the
equation. 
We   begin (in Section~\ref{subsec:setting}) by writing down in an
informal way the equations and the expected expansion 
 of the solution around a given state. That will motivate the
 computations of Section~\ref{subsec:lemma} in which an abstract
 result 
 is proved, showing under what assumptions on the equation one can
 justify such an expansion.

\subsection{Setting of the problem}
\label{subsec:setting}
Consider a first order partial differential equation, of the form
\begin{equation*}
  \d_t u = L(\d_x)u,\quad (t,x)\in \R\times\R^n,\quad u:\R\times\R^n\to
  \C\text{ or }\R^d, \ d\ge 1. 
\end{equation*}
We assume that the evolution of the solution to this linear equation
is described by a group $U(t)$. 
In the semilinear equations we have in mind, the nonlinearity will be
 a power law  $\Phi$
 of degree~$p \ge 2$. Let us
 consider   any  solution~$\ul $   to 
 the following equation
 \begin{equation*}
   \d_t \ul =L(\d_x)\ul + \Phi(\ul).
 \end{equation*}
Introduce the Duhamel formula associated to this equation:
 \begin{equation}
 \label{equl}
\ul (t) = U(t) \ul_{0} + N(\ul)(t),
 \end{equation}
 where we have defined
 \begin{equation}
   \label{eq:Duhamelinhom}
 N(\ul )(t) := \int_{t_0}^{t} U(t-s) \Phi\(\ul (s)\) \: ds.
 \end{equation}
 In scattering theory, one must think of the initial time as being
infinite, $t_0=-\infty$, in which case $\ul_{0}=u_-$ is an asymptotic
state. 
 \begin{example}\label{ex:Schrodinger}
   To make our discussion a little more concrete, we illustrate it
   with the case of a nonlinear Schr\"odinger equation
   \begin{equation}\label{eq:exNLS}
     i\d_t u + \frac{1}{2}\Delta u =\lvert u\rvert^{p-1}u. 
   \end{equation}
In this case, 
$U(t)=e^{i\frac{t}{2}\Delta}$, and $\Phi(u) = -i\lvert u\rvert^{p-1}u$. 
 \end{example}
 \begin{example}\label{ex:ondes}
   In the case of the nonlinear wave equation
   \begin{equation}\label{eq:exNLW}
     \d_t^2 u -\Delta  u + u^p=0,
   \end{equation}
we set $\ul = \,^t(u,\d_t u)$. Denote
\begin{equation*}
  \omega = \(-\Delta\)^{1/2}\ ;\ W(t)=\om^{-1}\sin \(\om t\) \ ;\ \dot
  W(t)=\cos\(\om t\).  
\end{equation*}
Then \eqref{eq:exNLW} takes the form
\eqref{equl}--\eqref{eq:Duhamelinhom}, with
\begin{equation*}
  U(t)=\(
  \begin{array}[c]{cc}
\dot W(t) & W(t)\\
-\om^2 W(t) & \dot W(t)
  \end{array}
\)\quad ;\quad \Phi(\ul)= \( 
\begin{array}[c]{c}
0\\
-\ul_1^p
\end{array}
\)= \( 
\begin{array}[c]{c}
0\\
-u^p
\end{array}
\).
\end{equation*}
The same holds in the case of the nonlinear Klein--Gordon equation
   \begin{equation*}
     \d_t^2 u -\Delta  u + u+u^p=0.
   \end{equation*}
The only adaptation needed in this case consists in substituting
$\omega$ with $\Lambda = (1-\Delta)^{1/2}$.  
 \end{example}
 We  suppose that this semilinear equation has global solutions in
 time and that a nonlinear scattering theory is available 
 (examples are provided in Section~\ref{sec:application} below). The
 discussion that follows is purely formal, and is intended 
 as a motivation to the computations carried out in the coming paragraph.
   
   \medskip
   
  Let us  
  construct a solution to the equation associated with an initial data
  which is   a perturbation of~$\ul_{0}$, 
   written~$\ul_0 + \eps u_{0}$  where~$\eps$ is a 
  small parameter, and   let us write
 the solution~$u^\eps$ under the form~$u^\eps = \ul + w^\eps$. We are
  looking for an 
  expansion of the perturbation~$w^\eps$ in powers  
 of~$\eps$. 
 Writing~$\Phi(\ul +w^\eps)$
 in terms of~$\Phi(\ul )$ using Taylor's formula yields easily that the
  equation on~$w^\eps$ must be of the following type: 
\begin{equation}
\label{eqformalw}
 w^\eps(t) = U(t) (\eps u_{0}) + \sum_{j = 1}^{p}\int_{t_0}^t U\(t-s\)\Phi_j\(
    \ul(s),w^\eps(s),\ldots,w^\eps(s)\) ds,
  \end{equation}
  where from now on~$\Phi_{j} (\alpha_{0},\alpha_{1},\dots,\alpha_{j})$
  denotes a multi-linear form,  
  which is~$(p-j)$-linear   in~$\alpha_{0}$ and linear
    in its~$j$ last arguments. In general, this multi-linearity is on
  $\R$ only, since in the case of the nonlinear Schr\"odinger
  equation, conjugation is involved in the above formula. To ease the
  notations, we introduce 
  \begin{equation}\label{eq:exprNj}
    N_{j} (\ul ,w,\dots,w)(t) = \int_{t_0}^t U\(t-s\)\Phi_j\(
    \ul(s),w(s),\ldots,w(s)\) ds.
  \end{equation}
  Our aim is now to write an expansion of~$w^\eps$ in powers of~$\eps$,
  $\displaystyle w^\eps = \sum_{k \in \N}\eps^{k+1}w_{k}$. 
  Two different situations can occur, according to the value
  of~$\ul_0$: either~$\ul_0$ is identically zero 
   (and the situation corresponds
  to the case of  small data), or it is not. 
  
  \subsubsection*{Case 1: Expansion around zero} Suppose~$\ul_0$
  vanishes identically. In that case the 
   only~$\Phi_{j}$   in~\eqref{eqformalw}
which is not identically zero  is when~$j = p$, and each~$w_{k}$ can
   be computed explicitly:  
 the only non vanishing terms in the expansion are   terms of the
   type~$w_{k(p-1)}$, for~$k \in \N$, with  
  $$
  w_{0}(t) =  U(t) u_{0}, 
  $$
  and where  the other terms of the expansion are given by an explicit
  algorithm, of the form
  \begin{equation*}
    w_{(k+1)(p-1)}(t)= G_k\(
    w_0(t),w_{p-1}(t),\ldots,w_{k(p-1)}(t)\),\quad k\ge 0. 
  \end{equation*}
 Typically, $w_0$ and $w_{p-1}$ are given by
 \begin{equation*}
   w_0(t) = U(t)u_0\quad ;\quad w_{p-1}(t)= N_p\(w_0(t),\ldots,w_0(t)\). 
 \end{equation*}
 \begin{example}
   In the above example of the nonlinear Schr\"odinger equation~\eqref{eq:exNLS}, this    yields
   \begin{align*}
     w_0(t,x)&= e^{i\frac{t}{2}\Delta}u_0(x),\\
 w_{p-1}(t,x)&=
     -i\int_{-\infty}^t e^{i\frac{t-s}{2}\Delta}\(\lvert
     e^{i\frac{s}{2}\Delta}u_0(x)
     \rvert^{p-1}e^{i\frac{s}{2}\Delta}u_0(x)\)ds. 
   \end{align*}
In other words, $w_0$ and $w_{p-1}$ solve
\begin{align*}
  &i\d_t w_0+\frac{1}{2}\Delta w_0=0\ &&;\quad
  U(-t)w_{0}(t)\big|_{t=-\infty}=u_0.\\ 
& i\d_t w_{p-1}+\frac{1}{2}\Delta w_{p-1}= \lvert
     w_0\rvert^{p-1}w_0\ &&;\quad U(-t)w_{1}(t)\big|_{t=-\infty}=0.
\end{align*}
It is obvious that $e^{-i\frac{t}{2}\Delta}w_0(t,x)$ converges as
$t\to +\infty$, and part of the game consists in showing  that so does
$e^{-i\frac{t}{2}\Delta}w_1(t,x)$.  
 \end{example}
  \subsubsection*{Case 2: Expansion around any initial data}
In that case all the~$\Phi_{j}$'s have to be taken into account
 in~\eqref{eqformalw}, so   the series will be 
 full if $\ul\not =0$. Moreover the~$w_{k}$'s are not
  computed explicitly. For instance the first two terms~$w_{0} $
 and~$w_{1}$ of the expansion satisfy 
 \begin{equation*}
   w_{0}(t) = U(t)u_{0}+ N_{1}\(\ul,w_{0}\)(t)\ \  ;\ \  w_{1}(t) =
   N_{1}\(\ul,w_{1}\)(t) + N_{2}\(\ul,w_{0},w_{0}\)(t). 
 \end{equation*}
 \begin{example}
   In our Schr\"odinger example \eqref{eq:exNLS}, this means that
   $w_0$ must solve 
   \begin{equation*}
\begin{aligned}
     i\d_t w_0 +\frac{1}{2}\Delta w_0 &= p \lvert \ul\rvert^{p-1} w_0 +
     (p-1) \ul^{ {(p+1)}/2}\overline \ul^{ (p-1 )/2}\overline w_0,\\
     U(-t)w_{0}(t)\big|_{t=-\infty}&=u_0.
\end{aligned}
   \end{equation*}
Note that the above Hamiltonian is not self-adjoint in
general. However, this aspect will not be an obstruction to our analysis.
 \end{example}
\subsubsection*{Conclusion} To summarize the above considerations, the
solution to the equation 
$$
u^\eps(t) = U(t) (\eps u_{0}) + N\(u^\eps\)(t)
$$
can be expanded as
$$
u^\eps =   \eps \sum_{k = 0}^{\infty} \eps^{k(p-1)} w_{k(p-1)},
$$
where   the~$w_{k(p-1)}$ satisfy  linear equations and can be computed
explicitly by induction. On the other hand, 
 the solution to the equation
$$
u^\eps(t) = U(t) (\ul_{0} + \eps u_{0}) + N\(u^\eps\)(t)
$$
can be expanded as
$$
u^\eps =  \ul +  \eps \sum_{k = 0}^{\infty} \eps^{k } w_{k },
$$
where again   the~$w_{k }$ satisfy  linear equations, but this time
are only known implicitly (again by induction). Those  
expansions allow to conclude that the scattering operator is analytic,
around any given state (small or large). 
In order to make those heuristical remarks rigorous, we need to prove
the convergence of the series formally obtained above. This 
is performed in the next section, where we prove an abstract result
stating under what conditions the series does converge.

\subsection{An abstract lemma}
\label{subsec:lemma}
In this section we adapt \cite[Theorem~3.2]{Masaki07} to
the case of a perturbation around any  
given state (in~\cite{Masaki07}, the perturbation is around zero only).

We keep the notation of the previous paragraph. Let us define~$D$ as
 the Banach  space in which  
 the    data lies, and~$F$ the space in which the linear flow
transports the   data.  The space~$F$ is a space-time Banach space,
 which we will write as~$F 
 = F_{1} \cap F_{2}$, where
$$F_{1}  :=  (C\cap L^\infty)(\R;D)$$
corresponds to the energy space, while
\begin{equation*}
 F_{2} =L^{q_1}(\R;X_1) \cap L^{q_2}(\R;X_2), \quad 1\le q_1,q_2<\infty.
\end{equation*}
for some Banach spaces~$X_1$ and $X_2$. Typically~$F_{2} $ should be
thought of as 
a Strichartz space, taking into account 
dispersive effects. In several applications, we will consider
$q_1=q_2$ and $X_1=X_2$. The main assumption on the linear evolution is that 
\begin{H1} There exists~$C_{0} >0$ such that for all~$ g \in D$,
\begin{equation*}
 \|U(\cdot)g\|_{F}
\le C_{0} \|g\|_{D}.  
\end{equation*}
\end{H1}
\noindent This assumption will always be satisfied thanks to Strichartz
estimates.   
\begin{example}
  Suppose that we consider the nonlinear Schr\"odinger at the $L^2$
  level. A natural choice is then $D=L^2(\R^n)$,
  $F= (C\cap L^\infty)(\R;L^2(\R^n))\cap L^q(\R;L^r(\R^n))$ for some
  Strichartz admissible pair $(q,r)$ (with $r=p+1$). 
\end{example}

As in the previous paragraph we consider a family of   $p$-linear 
forms denoted by~$(N_{j})_{1 \le j \le p}$,  
who are~$(p-j)$-linear in the first variable and linear in each of
the~$j$ remaining variables.  
We recall that the family~$(N_{j})_{1 \le j \le p}$ is constructed
as follows: 
\begin{equation}
\label{defNj}
\forall (a,b), \quad N(a+b) - N(b) =  \sum_{j = 1}^{p}N_{j} (a ,b,\dots,b).
\end{equation}
We will consider the second assumption:
\begin{H2}
There exists $\delta,C>0$ such that for all
$\ul,u_{1},\dots,u_{j}  \in F$ and for all $I$ interval in $\R$, we
have:
\begin{align*}
 & \|{\mathbf 1}_{t \in I} N_{j}(\ul ,u_{1},\dots, 
u_{ j})\|_{F} \le     C \|{\mathbf 1}_{t \in I}
\ul\|_{F_{2}}^\delta 
\|\ul\|_{F}^{p-\delta-j}
\prod_{\ell = 1}^{j} \|u_{\ell}\|_{F}\text{ if }j\le p-1,\\
 & \|{\mathbf 1}_{t \in I} N_{p}(u_{1},\dots, 
u_{ p})\|_{F} \le     C \sum_{\ell=1}^p\|{\mathbf 1}_{t \in I}
 u_\ell\|_{F_{2}}^\delta  \|
 u_\ell\|_{F}^{1-\delta} 
\prod_{\ell'\not= \ell}^{p} \|u_{\ell' }\|_{F}
.
\end{align*}
\end{H2}

\begin{remark}
The definition of~$F$ implies that
 if~$A$ and~$B$ are two disjoint intervals
of~$\R$, then
\begin{equation}\label{H2}
\|{\mathbf 1}_{t \in A \cup B} f\|_{F} \sim \|{\mathbf
  1}_{t \in A  } f\|_{F} + \|{\mathbf 1}_{t \in B  } f\|_{F}.  
\end{equation}
Moreover Lebesgue's theorem implies that
\begin{equation}\label{H6}
\forall v \in F_{2}, \quad \lim_{T \to  +\infty} \|{\mathbf 1}_{t \ge
  T} v\|_{F_{2}} = 0. 
\end{equation}
Similarly, we notice that~($H2$), applied to~$j = 1$,
implies  that~$\R$ may be decomposed into  a finite, 
disjoint union of~$K$ intervals~$(I_{k})_{1 \le k  
\le K}$  such that
\begin{equation}\label{H4} \quad  \forall  \ul,v\in F, \quad
  \|{\mathbf 1}_{t \in I_{k}} N_{1}(\ul,v) 
\|_{F} \le \frac 12 \|{\mathbf 1}_{t \in I_{k}} v
\|_{F}.
\end{equation}
\end{remark} 

\smallbreak

Fix~$u_{0}$   in~$D$. We construct by induction
a family~$(w_{k})_{k \in \N}$:
$$
w_{0}(t) = U(t)u_{0} + N_{1}(\ul ,w_{0})(t),
$$
$$
w_{m} = \sum_{j = 1}^{p} \sumetage{j + \ell_{1}+\dots+\ell_{j} =
  m+1}{\ell_{i} \ge 0} 
N_{j}(\ul ,w_{ \ell_{1}},\dots,
w_{ \ell_{j}}),
$$
with the convention that~$\displaystyle\sum_{\emptyset} = 0$. 
We have the following important lemma.
\begin{lemma}\label{lem:meta}
Let~$\ul \in F$ solve~\eqref{equl}  with initial data~$\ul_{0} \in D$,
and let~$u_{0}$ be a given function in~$D$, with $\|u_0\|_D\le M$.  
Assume~($H1$) and~($H2$) hold.
Then there exists $\eps_0=\eps_0\(\|\ul\|_F,M\)>0$ such that 
for~$0<\eps\le \eps_0$, the 
series~$\displaystyle \sum_{k \in \N} \eps^{k}w_{k}$ converges
normally in~$F$, and
\begin{equation*}
 u^\eps := \ul + \eps\displaystyle \sum_{k \in \N} \eps^{k}w_{k}
\quad\text{solves:}\quad u^\eps(t) = U(t)(\ul_{0} + \eps u_0) +
N\(u^\eps\)(t). 
\end{equation*}
\end{lemma}
\begin{remark}
Lemma \ref{lem:meta} implies in particular the real analyticity of the
wave operators as functions of~$D$, by considering the above result at
time $t=0$, since for $t_0=-\infty$, $u^\eps_{\mid t=0}=W_-\( \ul_0
+\eps u_0\)$.  
\end{remark}
\begin{proof}  
Let us start by finding a bound on~$w_{0}$  in~$F$. Inequality~\eqref{H4}
allows to write that 
\begin{eqnarray*}
\|{\mathbf 1}_{t \in I_{k}} w_{0}\|_{F} &\le& \|{\mathbf 1}_{t \in
  I_{k}} U(\cdot)u_0\|_{F} + 
\|{\mathbf 1}_{t \in I_{k}} N_{1}(\ul,w_{0})\|_{F} \\
 &\le& \|{\mathbf 1}_{t \in I_{k}} U(\cdot)u_0\|_{F} + \frac12 \|{\mathbf
  1}_{t \in I_{k}} w_{0}\|_{F}. 
\end{eqnarray*}
This implies directly, using~\eqref{H2}, that
$$
\|w_{0}\|_{F} \lesssim \|U(\cdot)u_0\|_{F}
$$
so by~($H1$) we infer that
\begin{equation}
\label{estimatew0}
\|w_{0}\|_{F} \lesssim C_{0} \|u_0\|_{D}.
\end{equation}
We prove by induction   that there exists $\Lambda
\ge 1$ such that for all $m\ge 1$,
\begin{equation*}
\|w_{m}\|_{F} \le \Lambda
^{m}. \tag{$\text{R}_{m}$} 
\end{equation*}
We notice that if that is the case, then the 
convergence of the series~$\displaystyle \sum_{k \in \N} \eps^{k}w_{k}$ in~$F$
is obvious as soon as~$\eps \Lambda <1$.

\smallbreak

Let us start by proving~($\text{R}_{1}$). 
 We have by definition
$$
w_{1} = N_{1}(\ul ,w_{1}) + N_{2}(\ul ,w_{0},w_{0}),
$$
and the same argument as in the case of~$w_{0}$ gives
\begin{eqnarray*}
\|{\mathbf 1}_{t \in I_{k}} w_{1}\|_{F} &\le&  
\|{\mathbf 1}_{t \in I_{k}} N_{1}(\ul,w_{1})\|_{F} + 
\|{\mathbf 1}_{t \in I_{k}} N_{2}(\ul ,w_{0},w_{0})\|_{F} \\
 &\le&  \frac12 \|{\mathbf 1}_{t \in I_{k}} w_{1}\|_{F} + 
\|{\mathbf 1}_{t \in I_{k}} N_{2}(\ul ,w_{0},w_{0})\|_{F} .
\end{eqnarray*}
By~\eqref{H2}, we infer that
$$
\|w_{1}\|_{F} \lesssim \| N_{2}(\ul ,w_{0},w_{0})\|_{F} .
$$
The continuity property~($H2$) then implies that
$$
\|w_{1}\|_{F} \lesssim C_{2}\|\ul\|_{F}^{p-2}  \|w_{0}\|_{F}^{2}
$$
so finally by~\eqref{estimatew0}
$$
\|w_{1}\|_{F} \lesssim C_{2}\|\ul\|_{F}^{p-2}(C_{0} \|u_0\|_{D})^{2}.
$$
So  we can choose~$\Lambda  \gtrsim 1+ C_{2}\|\ul\|_{F}^{p-2}(C_{0}
\|u_0\|_{D})^{2}  $ to get 
$$
\|w_{1}\|_{F} \le \Lambda.
$$ Now let us turn to the hierarchy of equations on~$w_{m}$, for~$m
\ge 2$. Supposing that~($\text{R}_{\ell}$)
holds for all~$1 \le \ell \le m-1$, let us prove~($\text{R}_{m}$). 
To simplify the notation 
we define
$$
\widetilde N (\ul,w_{0},\dots,w_{m-1}) := \sum_{j = 2}^{p}
\sumetage{j + \ell_{1}+\dots+\ell_{j} = m+1}{\ell_{i} \ge 0}
N_{j}(\ul ,w_{ \ell_{1}},\dots,
w_{ \ell_{j}}).
$$ 
The same argument as above yields
\begin{eqnarray*}
\|{\mathbf 1}_{t \in I_{k}} w_{m}\|_{F} &\le& 
\|{\mathbf 1}_{t \in I_{k}} N_{1}(\ul,w_{m})\|_{F} + \|{\mathbf 1}_{t
  \in I_{k}} 
\widetilde N (\ul,w_{0},\dots,w_{m-1})\|_{F} \\
 &\le&   \frac12 \|{\mathbf 1}_{t \in I_{k}} w_{m}\|_{F} + \|{\mathbf
   1}_{t \in I_{k}} 
\widetilde N (\ul,w_{0},\dots,w_{m-1})\|_{F}.
\end{eqnarray*}
Obviously this implies, using~\eqref{H2}, that
$$
\|w_{m}\|_{F} \lesssim \| \widetilde N (\ul,w_{0},\dots,w_{m-1})\|_{F}.
$$
By~($H2$) and defining~$C := \max_{1 \le j \le p} C_{j}$, we get that
\begin{eqnarray*}
\|w_{m}\|_{F} & \lesssim & C \sum_{j = 2}^{p} \|\ul\|_{F}^{p-j}
\sumetage{j + \ell_{1}+\dots+\ell_{j} = m+1}{\ell_{i} \ge 0}
\prod_{i = 1}^{j} \|w_{\ell_{i}}\|_{F} \\
& \lesssim & C \sum_{j = 2}^{p}  \|\ul\|_{F}^{p-j} \sumetage{j +
  \ell_{1}+\dots+\ell_{j} 
  = m+1}{\ell_{i} \ge 0} 
\prod_{i = 1}^{j} \Lambda ^{\ell_{i}} (C_{0} \|u_0\|_{D})^{\sharp\{i,
  \ell_{i} = 0\}} \\ 
& \lesssim & C  (1+C_{0} \|u_0\|_{D} +  \|\ul\|_{F})^p \sum_{j =
  2}^{p}\Lambda^{m+1-j} \\ 
& \lesssim & C  (1+C_{0} \|u_0\|_{D} +\|\ul\|_{F})^p \Lambda^{m-1}
\end{eqnarray*}
since~$\Lambda\ge 1$. To summarize, choosing
$$
\Lambda   \gtrsim 1+ C(1+C_{0}\|u_0\|_{D} +\|\ul\|_{F})^{p} 
$$
we  have $\|w_{m}\|_{F} \le \Lambda ^{m},$
and~$  (\text{R}_{m})$ is proved for all~$m \ge 1$. As remarked above, this
enables us to infer  that as soon as~$\eps$ is small enough,   
the series of general term~$\eps^{k}w_{k}$ is convergent.
\medskip

To conclude the proof of the lemma, let us prove that the solution of
\begin{equation}
\label{equ}
u^\eps(t) = U(t)(\ul_{0} + \eps u_0) + N\(u^\eps\)(t) 
\end{equation}
satisfies
$$\displaystyle u^\eps =\ul + \eps\displaystyle \sum_{k \in \N} \eps^{k}w_{k}.
$$
We show that the solution~$u^\eps$ of~\eqref{equ} satisfies
$$
\lim_{n\to \infty} \left\lVert u^\eps -\ul- \eps \sum_{ k= 0}^{n}
\eps^{k}w_{k}\right\rVert_{F} = 0, 
$$
by writing the equation satisfied by~$\widetilde w^\eps_n:= u^\eps -\ul-
\eps \sum_{ 
  k= 0}^{n} \eps^{k}w_{k}$. It is here that the  exact definition  
of the multi-linear operators~$N_{j}$ given in~\eqref{defNj} is
used. First, we know that for $\eps\Lambda <1$, the series $\sum
\eps^k w_k$ converges 
normally in $F$. Therefore,  $\widetilde w^\eps_n$ has a
limit in $F$ as $n\to \infty$, provided that $\eps$ is fixed such that
$\eps\Lambda <1$. On the other hand, by the 
definition of~$\widetilde w^\eps_n$, 
\begin{align*}
  \widetilde w^\eps_n &= 
N\(\ul + \eps \sum_{k = 0}^{n} \eps^{k} w_k + \widetilde w^\eps_n\) -
 N(\ul) \\
&- \eps  \sum_{k =0}^{n} \eps^{k} \sum_{j = 1}^{p}
\sumetage{j + \ell_{1}+\dots+\ell_{j} = k+1}{\ell_{i} \ge 0}
 N_{j}(\ul ,w_{ \ell_{1}},\dots,
w_{ \ell_{j}})\\
&= \sum_{j=1}^p N_j\( \ul,\eps \sum_{\ell_1=0}^n \eps^{\ell_1}
w_{\ell_1} + \widetilde w^\eps_n,\ldots, \eps \sum_{\ell_j=0}^n \eps^{\ell_j}
w_{\ell_j} + \widetilde w^\eps_n\) \\
&- \eps  \sum_{k =0}^{n} \eps^{k} \sum_{j = 1}^{p}
\sumetage{j + \ell_{1}+\dots+\ell_{j} = k+1}{\ell_{i} \ge 0}
 N_{j}(\ul ,w_{ \ell_{1}},\dots,
w_{ \ell_{j}})
\end{align*}
From the above estimates, we can write
\begin{align*}
 \sum_{j=1}^p &N_j\( \ul,\eps \sum_{\ell_1=0}^n \eps^{\ell_1}
w_{\ell_1} + \widetilde w^\eps_n,\ldots, \eps \sum_{\ell_j=0}^n \eps^{\ell_j}
w_{\ell_j} + \widetilde w^\eps_n\)= G_n\(\widetilde w^\eps_n\) \\
&\quad \quad + \sum_{j=1}^p N_j\( \ul,\eps \sum_{\ell_1=0}^n \eps^{\ell_1}
w_{\ell_1} ,\ldots, \eps \sum_{\ell_j=0}^n \eps^{\ell_j}
w_{\ell_j}\),
\end{align*}
where $G_n$ is such that we can decompose $\R$ as a finite, disjoint
union of intervals $J_q$, $1\le q\le Q$, independent of $n$, such that  
\begin{equation*}
  \lVert {\mathbf 1}_{t\in J_q} G_n\(\widetilde w^\eps_n\)\rVert_F \le
  \frac{1}{2}\lVert {\mathbf 1}_{t\in J_q} \widetilde
  w^\eps_n\rVert_F.  
\end{equation*}
We infer
\begin{align*}
  \widetilde w^\eps_n = G_n\(\widetilde w^\eps_n\) +
  \sum_{k=n+1}^{p-1+pn} \eps^{1+k} \sum_{j=1}^p \sumetage{j +
  \ell_{1}+\dots+\ell_{j} = k+1}{0\le \ell_{i} \le n} 
 N_{j}(\ul ,w_{ \ell_{1}},\dots,
w_{ \ell_{j}}). 
\end{align*}
Using \eqref{H2} and summing over the intervals $J_q$, we conclude
\begin{equation*}
  \lVert  \widetilde
  w^\eps_n\rVert_F =\O\(\(\eps\Lambda\)^{n+2}\).
\end{equation*}
Since $\eps\Lambda<1$ in order for all the above estimates to hold,
Lemma~\ref{lem:meta} follows from uniqueness for \eqref{equ} in $F$, which in
turn is a consequence of  ($H2$). 
\end{proof}

This result allows us to infer the analyticity of the scattering
operator, as shown in the following lemma. 
\begin{lemma}\label{lem:scattering}
  Let the assumptions ($H1$) and ($H2$) be satisfied. 
Assume furthermore that~$U(\cdot)$ is uniformly continuous in~$D$. Then
 $U(-t)\ul(t)$ converges  to  a limit $\ul_{+} $ in $D$ as $t\to +\infty$, and
 for   all $k\ge 0$, $U(-t)w_k(t)$ has a limit in~$D$, denoted by
  $w_k^+$. Moreover, for   $\eps$ sufficiently small, the series $\sum_{k\in \N
  }\eps^{k} w_k^+$   converges normally in $D$ and the function
\begin{equation*}
  u_+^\eps :=\ul_+ + \eps\sum_{k\in \N} \eps^{k}w_k^+
\end{equation*}
is the limit of $U(-t)u^\eps(t)$ in $D$ as $t\to +\infty$. 

In particular, the scattering operator is  analytic from $D$ to $D$.
\end{lemma}
\begin{proof}
Let us start by proving the existence of~$\ul_{+}$. We have
\begin{eqnarray*}
\left\|U(-t_{2}) \ul (t_{2})- U(-t_{1}) \ul(t_{1})\right\|_{D} & = & \left\|
\int_{t_{1}}^{t_{2}} U(-s) \Phi(\ul) \: ds
\right\|_{D}
\\
& = & \left\|{\mathbf 1}_{[t_{1},t_{2}]}
U(-t)N(\ul)
\right\|_{D} \\
&\lesssim & \left\|{\mathbf 1}_{[t_{1},t_{2}]} N(\ul)
\right\|_{F_{1}} \le  \left\|{\mathbf 1}_{[t_{1},t_{2}]} N(\ul)
\right\|_{ F}\\
&\lesssim & \left\|{\mathbf 1}_{[t_{1},t_{2}]} \ul\right\|_{
  F_{2}}^\delta
\left\|\ul\right\|^{p-\delta}_{ F},
\end{eqnarray*}
by assumption~($H2$). We conclude by the fact that the right-hand side
goes to zero as~$t_{1},t_{2}$ go to infinity.   

Now we prove the result on~$U(-t)w_k(t)$ by induction on~$k$. For~$k =
0$ we have, in the same fashion as above,
\begin{eqnarray*}
\left\|U(-t_{2}) w_{0} (t_{2})- U(-t_{1}) w_{0}(t_{1})\right\|_{D} 
& = & \left\|{\mathbf 1}_{[t_{1},t_{2}]}
U(-t)N_{1}(\ul,w_{0})
\right\|_{D} \\
&\lesssim &   \left\|{\mathbf 1}_{[t_{1},t_{2}]} N_{1}(\ul,w_{0})
\right\|_{ F}\\
&\lesssim & \left\|{\mathbf 1}_{[t_{1},t_{2}]} \ul\right\|_{
  F_{2}}^\delta \left\| \ul\right\|^{p-1-\delta}_{
  F} \|w_{0}\|_{F}, 
\end{eqnarray*}
since~$w_{0}$ belongs to~$F$ due to~\eqref{estimatew0}. We conclude as above.

Now suppose that for~$m \ge 1$ and for all~$0 \le  \ell \le m-1 $, 
$U(-t)w_{\ell}(t)$ has a limit. We   prove the result for~$U(-t)w_m(t)$.
We have as above
\begin{align*}
 \|U(-t_{2}) w_{m} (t_{2}) &-  U(-t_{1}) w_{m}(t_{1})\|_{D} 
  \le  \\
& \le  \sum_{j = 1}^{p} \sumetage{j + \ell_{1}+\dots+\ell_{j} =
  m+1}{\ell_{i} \ge 0} 
\left\|{\mathbf 1}_{[t_{1},t_{2}]}
U(-t)N_{j}(\ul ,w_{ \ell_{1}},\dots,
w_{ \ell_{j}})
\right\|_{D} \\
&\lesssim   \sum_{j = 1}^{p} \sumetage{j + \ell_{1}+\dots+\ell_{j} =
  m+1}{\ell_{i} \ge 0} 
  \left\|{\mathbf 1}_{[t_{1},t_{2}]}N_{j}(\ul ,w_{ \ell_{1}},\dots,
w_{ \ell_{j}})
\right\|_{ F}\\
&\lesssim   \sum_{j = 1}^{p-1} \sumetage{j + \ell_{1}+\dots+\ell_{j} =
  m+1}{\ell_{i} \ge 0} 
 \left\|{\mathbf 1}_{[t_{1},t_{2}]} \ul\right\|^{\delta}_{
  F_{2}}\left\|\ul\right\|^{p-\delta-j}_{
  F}\prod_{k = 1}^{j} 
 \|w_{\ell_{k}}\|_{F}\\
&+ \sumetage{p + \ell_{1}+\dots+\ell_{p} =
  m+1}{\ell_{i} \ge 0} \sum_{k=1}^p
 \left\|{\mathbf 1}_{[t_{1},t_{2}]} u_{\ell_k}\right\|^{\delta}_{
  F_{2}}\left\|u_{\ell_k}\right\|^{1-\delta}_{
  F}\prod_{k'\not = k}
 \|w_{\ell_{k'}}\|_{F}.
\end{align*}
The result follows as previously.

The convergence of the series defining~$u^{\eps}_{+}$ is due to
Lemma~\ref{lem:meta}, and  
Lemma~\ref{lem:scattering} follows directly.
\end{proof}

\subsection{An easy and useful adaptation}
\label{subsec:adaptlemma}

For nonlinear Schr\"odinger and wave equations, Lemmas~\ref{lem:meta}
and \ref{lem:scattering} are well adapted to 
study the wave and scattering operators in energy spaces. On the other
hand, as recalled in the introduction, weighted Sobolev spaces are
very useful in scattering theory for these equations. Typically, for
the nonlinear Schr\"odinger equation, the natural energy space is
$H^1(\R^n)$, but more results concerning scattering are available in
$\Sigma$, defined in the introduction. In the case of the energy space
$H^1$, we will see that the natural choice for the space $F$ is
\begin{equation*}
  F = C\cap L^\infty \( \R;H^1(\R^n)\)\cap L^{\frac{4p+4}{n(p-1)}}\(\R;
  W^{1,p+1}(\R^n)\), 
\end{equation*}
which is of the form considered in \S\ref{subsec:lemma}, with $X=
W^{1,p+1}(\R^n)$. When working on $\Sigma$, the natural choice for $F$
is
\begin{equation*}
 \widetilde F = F\cap \left\{ f\in C(\R;\Sigma), \ J(t)f \in
  L^{\frac{4p+4}{n(p-1)}}\(\R; L^{p+1}(\R^n)\)\right\},
\end{equation*}
where $J(t)= x+it\nabla$ is the Galilean operator. It satisfies the
important property $J(t)= U(t)xU(-t)$. The situation is
fairly similar in the case of the nonlinear wave equation. 
\smallbreak

It is therefore natural to adapt the framework of
\S\ref{subsec:lemma}. For the same spaces~$D$ and~$F$, introduce
\begin{equation*}
  \widetilde F = F\cap F_3,\quad \text{where }\lVert f\rVert_{F_3} =
 \lVert J f\rVert_{L^\infty(\R;E)} +\lVert J
 f\rVert_{L^q(\R;Y)}, 
\end{equation*}
for some Banach spaces $E$ and $Y$, and some operator $J$
depending on time. Define the space
$\widetilde D$ and $\widetilde F_2$ by their norms
\begin{equation*}
 \lVert g\rVert_{\widetilde D} = \lVert g\rVert_{D}+
\lVert J(0)g\rVert_{E} \quad ;\quad
\lVert f\rVert_{\widetilde F_2} = \lVert f\rVert_{F_2}+
\lVert J f\rVert_{L^q(\R;Y)}. 
\end{equation*}
It is easy to check that Lemmas~\ref{lem:meta}
and \ref{lem:scattering} remain valid if $F$ is replaced
by~$\widetilde F$, provided that  (H1)  and  (H2)  
 are replaced by:
\begin{equation*}
\exists C_{0} , \quad \forall g \in \widetilde D,  \quad 
\|U(\cdot)g\|_{\widetilde F} 
\le C_{0} \|g\|_{\widetilde D}. \tag{ $\widetilde{ H 1}$} 
\end{equation*}
and
\begin{H2tilde}
  There exists $\delta,C>0$ such that for all
$\ul,u_{1},\dots,u_{j}  \in \widetilde F$ and for all $I$ interval in $\R$, we
have:
\begin{align*}
 & \|{\mathbf 1}_{t \in I} N_{j}(\ul ,u_{1},\dots, 
u_{ j})\|_{\widetilde F} \le     C \|{\mathbf 1}_{t \in I}
\ul\|_{\widetilde F_{2}}^\delta 
\|\ul\|_{\widetilde F}^{p-\delta-j}
\prod_{\ell = 1}^{j} \|u_{\ell}\|_{\widetilde F}\text{ if }j\le p-1,\\
 & \|{\mathbf 1}_{t \in I} N_{p}(u_{1},\dots, 
u_{ p})\|_{\widetilde F} \le     C \sum_{\ell=1}^p\|{\mathbf 1}_{t \in I}
 u_\ell\|_{\widetilde F_{2}}^\delta  \|
 u_\ell\|_{\widetilde F}^{1-\delta} 
\prod_{\ell'\not= \ell}^{p} \|u_{\ell' }\|_{\widetilde F}
.
\end{align*}
\end{H2tilde}
In the applications, we shall also use the following lemma, whose
proof follows the same lines as the proofs of Lemmas~\ref{lem:meta}
and \ref{lem:scattering}, and is left out.
\begin{lemma}\label{lem:generalized}
  Let~$\ul \in \widetilde F$ solve~\eqref{equl}  with initial
  data~$\ul_{0} \in \widetilde D$, 
and let~$u_{0}$ be a given function in~$\widetilde D$, with
  $\|u_0\|_{\widetilde  D}\le M$.  
Assume~($\widetilde {H 1}$) and~($\widetilde{H 2}$) hold.
Then there exists $\eps_0=\eps_0\(\|\ul\|_F,M\)>0$ such that 
for~$0<\eps\le \eps_0$, the 
series~$\displaystyle \sum_{k \in \N} \eps^{k}w_{k}$ converges
normally in~$\widetilde F$, and
\begin{equation*}
 u^\eps := \ul + \eps\displaystyle \sum_{k \in \N} \eps^{k}w_{k}
\quad\text{solves:}\quad u^\eps(t) = U(t)(\ul_{0} + \eps u_0) +
N\(u^\eps\)(t). 
\end{equation*}
Assume furthermore that~$U(\cdot)$ is uniformly continuous
 in~$\widetilde D$. Then
 $U(-t)\ul(t)$ converges  to  a limit $\ul_{+} $ in $\widetilde D$ as
 $t\to +\infty$, and 
 for   all $k\ge 0$, $U(-t)w_k(t)$ converges to $w_k^+$ in~$\widetilde
 D$. Moreover, for   $\eps$ sufficiently small, the series $\sum_{k\in
 \N 
  }\eps^{k} w_k^+$   converges in $\widetilde D$ and the function
\begin{equation*}
  u_+^\eps :=\ul_+ + \eps\sum_{k\in \N} \eps^{k}w_k^+
\end{equation*}
is the limit of $U(-t)u^\eps(t)$ in $\widetilde D$ as $t\to +\infty$. 

In particular, the scattering operator is analytic 
from~$\widetilde D$ to $\widetilde D$. 
\end{lemma}
\section{Application to semilinear dispersive equations}
\label{sec:application}

\subsection{The Schr\"odinger equation}
\label{subsec:NLS}

\subsubsection{General presentation}

We consider the nonlinear Schr\"odinger equation with gauge invariant
nonlinearity presented in the introduction:
\begin{equation}
  \label{eq:NLS}
  i\d_t u+\frac{1}{2}\Delta u = \lvert u\rvert^{p-1}u,\quad (t,x)\in
  \R\times \R^n.
\end{equation}
In order for the nonlinearity to be analytic, we assume that $p$ is an
odd integer, with $p\ge 3$. Note that compared to
Eq.~\eqref{eq:NLSintro}, we have imposed the value $\lambda=+1$ for
the coupling constant. We consider defocusing nonlinearities, for
which the scattering theory is much richer than in the focusing case,
where the existence of solitons and finite time blow-up phenomenon
may prevent the solution $u$ from scattering at infinity. 
\smallbreak

Two different frameworks seem particularly well suited to study
scattering for \eqref{eq:NLS}: $H^1(\R^n)$, and 
\begin{equation*}
  \Sigma = \{ f\in H^1(\R^n),\quad
  x\mapsto \lvert x\rvert f(x)\in L^2(\R^n)\}.
\end{equation*}
We apply Lemmas~\ref{lem:meta} and \ref{lem:scattering} in the first
case, and Lemma~\ref{lem:generalized} in the second case. Note that
another framework should be well suited as well, which is the $L^2$
case. If $p>1+4/n$, then the nonlinearity in \eqref{eq:NLS} is
$L^2$-supercritical: the results of \cite{CCT2} show that a scattering 
theory in $L^2$ with continuous dependence on the data is hopeless. If
$p<1+4/n$, then scattering is not known at the $L^2$ level, and does
fail if $p\le 1+2/n$ (\cite{Barab,Ginibre,Strauss74,Strauss81}). In
the $L^2$-critical case $p=1+4/n$, scattering is known for small data
\cite{CW89}. Note that $p=1+4/n$ is an odd integer only when $n=1$ or
$2$. For $n=1$, scattering for large $L^2$ data is not known so
far. For $n=2$, scattering for large $L^2$ \emph{radial} data was
proved in \cite{KTV07}. To avoid an endless numerology, we leave out
the discussion on the $L^2$ case at this stage.  
\smallbreak

Note also that the case of non-Euclidean geometries could be
considered. In \cite{BCS}, the existence of scattering
operators was established in $H^1$ for solutions to the nonlinear
Schr\"odinger equation on hyperbolic space, in space dimension three,
for  
energy-subcritical nonlinearities: the nonlinearity is analytic if it
is cubic (and only in that case, since the energy-critical case has
not been treated so far). Also, from the results in \cite{IS08},
scattering in $H^1$ is available on the two-dimensional hyperbolic
space. The 
analyticity of wave and scattering operators in these cases can then
be established 
by  the same argument as in \S\ref{sec:NLSH1} below. 
\subsubsection{The case of $H^1$}
\label{sec:NLSH1}
For $p\ge 1+4/n$, with $p<1+4/(n-2)$ when $n\ge 3$, the existence and
continuity of wave operators was established in \cite{GV85}. If we assume
moreover that $p>1+4/n$, then asymptotic completeness holds: this was
proved initially in \cite{GV85}
for $n\ge 3$ (see also \cite{TaoVisanZhang} for a simplified proof), and in
\cite{Nakanishi99,Nakanishi01} for $n=1,2$ (see also 
\cite{CazCourant}). 
We assume $1+4/n<p< 1+4/(n-2)$. In
order to prove the second part of Theorem~\ref{theo:NLS} in the
energy-subcritical case, it suffices 
to exhibit  spaces $D$ and $F_2$ such that ($H1$) and ($H2$) are
satisfied. We consider the energy-critical case $p=1+4/(n-2)$ in a
different paragraph, since the proof is slightly different.
\smallbreak

We set naturally $D=H^1(\R^n)$, hence $F_1=(C\cap
L^\infty)(\R;H^1(\R^n))$. 
The space $F_2$ is motivated by Strichartz estimates:
\begin{equation*}
  F_2 = L^{\frac{4p+4}{n(p-1)}}\(\R;W^{1,p+1}(\R^n)\).
\end{equation*}
Note that the pair $(q,r)=(\frac{4p+4}{n(p-1)}, p+1)$ is
$L^2$-admissible:
\begin{equation*}
  \frac{2}{q}=n\(\frac{1}{2}-\frac{1}{r}\)=:\delta(r),\quad 2\le
  r\le \frac{2n}{n-2}, \ (n,q,r)\not = (2,2,\infty). 
\end{equation*}
The fact that ($H1$) is satisfied is a consequence of homogeneous Strichartz
inequalities (\cite{GV85b,KT}). To check ($H2$), 
we use inhomogeneous Strichartz
inequalities, and the following
algebraic lemma:
\begin{lemma}\label{lem:alg1}
  Let $p\ge 1+4/n$, with $p< 1+4/(n-2)$ if $n\ge 3$. Set
  \begin{equation*}
    (q,r)=\(\frac{4p+4}{n(p-1)}, p+1\).
  \end{equation*}
Then $(q,r)$ is admissible. Set
$$\theta =\frac{p+1}{p-1}\times\frac{n(p-1)-4}{n(p-1)} .$$
Then $\theta \in [0,1[$. Define $s=r=p+1$ and $k=q/(1-\theta)$. Obviously,
\begin{equation*}
\frac{1}{s}=\frac{1-\theta}{r}+
\frac{\theta}{p+1}\quad ;\quad  
    \frac{1}{k}=\frac{1-\theta}{q}+
\frac{\theta}{\infty}\virgp
\end{equation*}
and we have: $\displaystyle  \frac{1}{r'}=\frac{1}{r}+\frac{p-1}{s}$,
  and $\displaystyle \frac{1}{q'}=\frac{1}{q}+\frac{p-1}{k} $. 
\end{lemma}
Recall that the nonlinear terms $N_j$ stem from an
inhomogeneous term in integral form, \eqref{eq:exprNj}. For a time
interval $I\subset \R$, inhomogeneous Strichartz estimates yield, for
$1\le j\le p$, 
\begin{equation*}
  \left\lVert {\mathbf 1}_{t\in I}N_j\( \ul, u_1,\ldots
  ,u_j\)\right\rVert_{L^\infty (\R;L^2)\cap L^q(\R;L^r)} \le C
  \left\lVert {\mathbf 1}_{t\in I}\left\lvert  \ul \right\rvert^{p-j}
  \prod_{\ell =1}^j \left\lvert u_\ell
  \right\rvert \right\rVert_{L^{q'}(\R;L^{r'})},
\end{equation*}
for some constant $C$ independent of $I$, and $\ul,u_1,\ldots,u_j\in
F$. Using Lemma~\ref{lem:alg1}, we infer, if $j\le p-1$:
\begin{align*}
  \left\lVert {\mathbf 1}_{t\in I}N_j\( \ldots
  \)\right\rVert_{L^\infty L^2\cap L^q L^r}& \lesssim
  \left\lVert {\mathbf 1}_{t\in I}\ul \right\rVert_{L^{q}L^{r}}
\left\lVert {\mathbf 1}_{t\in I}\ul
  \right\rVert_{L^{k}L^{s}}^{p-1-j} \prod_{\ell=1}^j \left\lVert
  {\mathbf 1}_{t\in I}u_\ell \right\rVert_{L^{k}L^{s}} \\
& \lesssim
  \left\lVert {\mathbf 1}_{t\in I}\ul \right\rVert_{L^{q}L^{r}} 
\left\lVert \ul
  \right\rVert_{L^{q}L^{r}}^{(1-\theta)(p-1-j)} \left\lVert \ul 
  \right\rVert_{L^{\infty}L^{p+1}}^{\theta(p-1-j)}\times\\
&\quad \times  \prod_{\ell=1}^j \left\lVert
  u_\ell \right\rVert_{L^{q}L^{r}}^{1-\theta}\left\lVert
  u_\ell \right\rVert_{L^{\infty}L^{p+1}}^{\theta}
\end{align*}
Using the embedding $H^1(\R^n)\hookrightarrow L^{p+1}(\R^n)$, we
deduce:
\begin{align*}
  \left\lVert {\mathbf 1}_{t\in I}N_j\( \ldots
  \)\right\rVert_{L^\infty L^2\cap L^q L^r}& \lesssim \left\lVert
  {\mathbf 1}_{t\in I}\ul \right\rVert_{L^{q}L^{r}} \left\lVert
  \ul \right\rVert_{F}^{p-1-j}   \prod_{\ell=1}^j \left\lVert 
u_\ell \right\rVert_{F}. 
\end{align*}
The estimate for $N_p$ in $L^\infty L^2\cap L^q L^r$ follows from the
same computation. To estimate $N_j$ in $L^\infty H^1\cap L^q W^{1,r}$,
we mimic the above computation. To simplify the presentation, and to
explain why Assumption~($H2$) is stated in such an apparently
intricate way, we consider only the case $j=1$. All the other cases
can be deduced in the same fashion. We have obviously
\begin{equation*}
  \left\lvert {\mathbf 1}_{t\in I}\nabla N_1\(\ul,u_1\)\right\rvert
  \lesssim \left\lvert {\mathbf 1}_{t\in I}
  \ul^{p-2}u_1\nabla \ul \right\rvert + \left\lvert {\mathbf
  1}_{t\in I}
  \ul^{p-1}\nabla u_1 \right\rvert.
\end{equation*}
Proceeding as above, we consider the $L^\infty L^2\cap L^q L^r$ norm,
and use H\"older's inequality, as suggested by
Lemma~\ref{lem:alg1}. However, we do not have the same room 
to balance the different Lebesgue's norms: we do not want to use
Sobolev embedding to control the derivatives. We find
\begin{align*}
  \left\lVert {\mathbf 1}_{t\in I}\nabla
  N_1\(\ul,u_1\)\right\rVert_{L^\infty L^2\cap L^q L^r}& \lesssim  
  \left\lVert \nabla \ul \right\rVert_{L^q L^r}
  \left\lVert {\mathbf 1}_{t\in I}\ul \right\rVert_{L^kL^s}^{p-2} 
 \left\lVert {\mathbf 1}_{t\in I}u_1 \right\rVert_{L^kL^s}\\
\ & +
\left\lVert \nabla u_1 \right\rVert_{L^q L^r}
  \left\lVert {\mathbf 1}_{t\in I}\ul \right\rVert_{L^kL^s}^{p-1}\\
 & \lesssim  
  \left\lVert \ul  \right\rVert_{F}
  \left\lVert {\mathbf 1}_{t\in I}\ul
  \right\rVert_{L^qL^r}^{(1-\theta)(p-2)}
 \left\lVert \ul
  \right\rVert_{L^\infty L^{p+1}}^{\theta(p-2)}
 \left\lVert u_1 \right\rVert_{F}\\
\ & +
\left\lVert u_1 \right\rVert_{F}
  \left\lVert {\mathbf 1}_{t\in I}\ul
  \right\rVert_{L^qL^r}^{(1-\theta)(p-1)}
\left\lVert \ul
  \right\rVert_{L^\infty L^{p+1}}^{\theta(p-1)}\\
&\lesssim \left\lVert {\mathbf 1}_{t\in I}\ul
  \right\rVert_{L^qL^r}^{1-\theta}\left\lVert \ul
  \right\rVert_{F}^{p+\theta -2} \left\lVert u_1 \right\rVert_{F},
\end{align*}
where we have used the same estimates as above (recall that $p\ge
3$). Therefore, Assumption~($H2$) is satisfied, with $\delta =
1-\theta$. Note that $\delta>0$ because we consider the
energy-subcritical case, $p<1+4/(n-2)$.  
\smallbreak

Therefore, we can apply Lemmas~\ref{lem:meta} and \ref{lem:scattering}
with $F$ as above. This yields the second part of
Theorem~\ref{theo:NLS}, except for the energy-critical case. Note that
in the following two cases:
\begin{itemize}
\item $n=1$ and $p=5$ (quintic nonlinearity),
\item $n=2$ and $p=3$ (cubic nonlinearity),
\end{itemize}
which are $L^2$ critical $p=1+4/n$, Lemma~\ref{lem:meta} shows that
the wave operators are analytic on $H^1(\R^n)$. However, scattering in
the energy space for arbitrary data is not known in these cases.

\subsubsection{The case of $\Sigma$}\label{sec:Sigma}

To overcome the drawback mentioned at the end of the previous
paragraph, we shall consider the weighted Sobolev space
$\Sigma$. Generally speaking, working in $\Sigma$ makes it possible to
decrease the admissible values for $p$ in order to have scattering,
from $p>1+4/n$, to $p\ge p_0(n)$, for some $1+2/n<p_0(n)<1+4/n$; see
\cite{CW92,GV79Scatt,HT87,NakanishiOzawa}. However, the gain in the
present context is rather weak, since we consider only integer values
for $p$: the gain corresponds exactly to the two cases pointed out
above. 
\smallbreak

As suggested in \S\ref{subsec:adaptlemma}, we consider the space 
\begin{equation*}
 \widetilde F = F\cap \left\{ f\in C(\R;\Sigma), \ J(t)f \in
  L^{\frac{4p+4}{n(p-1)}}\(\R; L^{p+1}(\R^n)\)\right\},
\end{equation*}
where $J(t)=x+it\nabla$, and $F$ was defined in the previous
paragraph. We can then mimic the 
above computation, in order to apply
Lemma~\ref{lem:generalized}. We recall two important properties of the
operator $J$ which make it possible to check
Assumptions~($\widetilde{H1}$) and ($\widetilde{H2}$): 
\begin{itemize}
\item It commutes with the linear Schr\"odinger group: $J(t)= U(t)x
  U(-t)$.
\item It acts on gauge invariant nonlinearities like a derivative,
  since
  \begin{equation*}
    J(t) = it e^{i\lvert x\rvert^2/(2t)}\nabla\(
  e^{-i\lvert x\rvert^2/(2t)} \cdot\),\quad \forall t\not =0.
  \end{equation*}
\end{itemize}
Lemma~\ref{lem:generalized} and the results of \cite{GV79Scatt} yield
Theorem~\ref{theo:NLS} in all the cases, but the energy critical one,
which is considered in the next paragraph. 

\subsubsection{The energy-critical case}

To complete the proof of Theorem~\ref{theo:NLS}, two cases  remain,
which correspond to the case $p=1+4/(n-2)$:
\begin{itemize}
\item $n=3$ and $p=5$.
\item $n=4$ and $p=3$.
\end{itemize}
Global existence and scattering for arbitrary data in $H^1(\R^n)$
were established in \cite{CKSTTAnnals} and \cite{RV}, respectively. 
A crucial tool in the energy critical case is the existence of
Strichartz estimates for $\dot H^1$-admissible pairs, as opposed to
the notion of $L^2$-admissible pairs used above. It is fairly natural
that our definition for $F$ is adapted in view of this notion. Recall
that for $n\ge 3$, a pair $(q,r)$ is $\dot H^1$-admissible if
\begin{equation*}
  \frac{2}{q}+\frac{n}{r}= \frac{n}{2}-1. 
\end{equation*}
Denote 
\begin{equation*}
  \gamma_0 = 2+\frac{4}{n} \quad \text{and}\quad
  \gamma_1=2+\frac{8}{n-2}. 
\end{equation*}
The pair $(\gamma_0,\gamma_0)$ is $L^2$-admissible, and
 $(\gamma_1,\gamma_1)$ is $\dot H^1$-admissible. We set
 \begin{align*}
   F=F_1\cap F_2,\text{ with }F_1& =(C\cap L^\infty)\(\R;H^1(\R^n)\),
   \text{ and}\\
F_2&=L^{\gamma_0}\(\R;W^{1,\gamma_0}(\R^n)\) \cap
   L^{\gamma_1}(\R\times\R^n).  
 \end{align*}
With such a space $F$, Assumption~($H1$) is satisfied, thanks to
Strichartz estimates, along with the Sobolev embedding $\dot
H^1(\R^n)\hookrightarrow L^{2n/(n-2)}(\R^n)$. To check that
Assumption~($H2$) is satisfied as well, we distinguish the two cases
we consider, for a more convenient numerology.

\subsubsection*{The quintic case, with $n=3$} In this case, we have
$\gamma_0=10/3$ and $\gamma_1=10$. For $u_1,\ldots,u_5\in F$, we have,
for $k=0$ or $1$, thanks to Strichartz estimates and H\"older's inequality:
\begin{align*}
  \Big\rVert \nabla^k\int_{t_0}^t & U(t-s)\(u_1\times\ldots\times
  u_5\)(s)ds\Big\rVert_{L^\infty(I;L^2)\cap L^{10/3}(I\times \R^n)}
  \\
& \lesssim \left\rVert \nabla^k\(u_1\times\ldots\times
  u_5\)\right\rVert_{L^{10/7}(I\times \R^n)}\\
&  \lesssim \sum_{j=1}^5 \left\rVert \nabla^k u_j\right\rVert_{L^{10/3}(I\times
  \R^n)}\prod_{\ell\not =j} \left\rVert u_\ell\right\rVert_{L^{10}(I\times
  \R^n)}
\lesssim \prod_{j=1}^5 \left\rVert
{\mathbf 1}_{t\in I}  u_j\right\rVert_{F_2}. 
\end{align*}
We also have, in view of Sobolev embedding,
\begin{align*}
  &\Big\rVert\int_{t_0}^t U(t-s)\(u_1\times\ldots\times
  u_5\)(s)ds\Big\rVert_{L^{10}(I\times \R^n)}
  \\
\lesssim &\Big\rVert\int_{t_0}^t U(t-s)\(u_1\times\ldots\times
  u_5\)(s)ds\Big\rVert_{L^{10}(I;W^{1,30/13})}\\
\lesssim &\sum_{k=0}^1\left\rVert \nabla^k\(u_1\times\ldots\times
  u_5\)\right\rVert_{L^{10/7}(I\times \R^n)},
\end{align*}
thanks to Strichartz estimates. Using the above computation, we infer
that Assumption~($H2$) is satisfied. 
 
\subsubsection*{The cubic case, with $n=4$}
In this case we have
$\gamma_0=3$ and $\gamma_1=6$. For $u_1,u_2,u_3\in F$, we have,
for $k=0$ or $1$, thanks to Strichartz estimates and H\"older's inequality:
\begin{align*}
  \Big\rVert {\mathbf 1}_{t\in I} \nabla^k\int_{t_0}^t &U(t-s)\(u_1 u_2
  u_3\)(s)ds\Big\rVert_{L^\infty_t L^2_x \cap L^{3}_{t,x}} 
   \lesssim \left\rVert \nabla^k\(u_1 u_2 u_3
  \)\right\rVert_{L^{3/2}(I\times \R^n)}\\
&  \lesssim \sum_{j=1}^3 \left\rVert \nabla^k u_j\right\rVert_{L^{3}(I\times
  \R^n)}\prod_{\ell\not =j} \left\rVert u_\ell\right\rVert_{L^6(I\times
  \R^n)}
\lesssim \prod_{j=1}^3 \left\rVert
{\mathbf 1}_{t\in I}  u_j\right\rVert_{F_2}. 
\end{align*}
We also have, in view of Sobolev embedding,
\begin{align*}
  \Big\rVert\int_{t_0}^t U(t-s)\(u_1 u_2 u_3
  \)&(s)ds\Big\rVert_{L^{6}(I\times \R^n)} 
  \\
\lesssim \Big\rVert\int_{t_0}^t U(t-s)\(u_1 u_2
  u_3\)&(s)ds\Big\rVert_{L^{6}(I;W^{1,12/5})}
\lesssim \sum_{k=0}^1\left\rVert \nabla^k\(u_1 u_2
  u_3\)\right\rVert_{L^{3}(I\times \R^n)}, 
\end{align*}
thanks to Strichartz estimates. Using the above computation, we infer
that Assumption~($H2$) is satisfied. 

Finally, it is easily checked that we can replace $H^1$ with $\Sigma$,
as in the previous paragraph. This completes the proof of
Theorem~\ref{theo:NLS}. 

\begin{remark}\label{rem:NLSBesov}
  At the level of $H^1$, it is possible to have a unified
  presentation, that is, without distinguishing the $H^1$-subcritical
  and $H^1$-critical cases. The price to pay consists in considering  
  Besov spaces for the definition of $F_2$, instead of Sobolev
  spaces. We have chosen to work in Sobolev for the simplicity and the
  explicit form of the computations. A more synthetic approach would
  consist in setting
\begin{align*}
  &F_2 = L^{\g_0} \(\R;B^{1}_{\g_0,2}(\R^n)\)\cap L^{\g_1}\(\R\times
  \R^n\),\\
\text{with }& \g_0=2+\frac{4}{n}\text{ and
  }\frac{p-1}{\g_1}+\frac{1}{\g_0} = \frac{1}{\g_0'}. 
\end{align*}
Sobolev and Strichartz inequalities are
replaced by
\begin{equation*}
  \lVert u \rVert_{L^\infty (\R;H^1)} +  \lVert u \rVert_{F_2} \le C\(
  \lVert u_0 \rVert_{H^1} + \Big\lVert i\d_t u +\frac{1}{2}\Delta
  u\Big\rVert_{L^{\g_0'}\big(\R; B^{1}_{\g_0',2}(\R^n)\big)}\),
\end{equation*}
an estimate established in \cite[\S3]{Nakanishi99}. Note that in the
energy-critical case $p=1+\frac{4}{n-2}$, this is  the estimate
which we have used, up the replacing Besov spaces $B_{p,2}^1$ with
$W^{1,p}$ (a modification which is non-trivial since $p\not =2$). 
\end{remark}

\subsection{The Hartree equation}
\label{subsec:hartree}

We now consider the Hartree equation \eqref{eq:r3intro} with a
defocusing nonlinearity, $\lambda =+1$, in space dimension $n\ge 3$: 
\begin{equation}
  \label{eq:r3}
  i\d_t u +\frac{1}{2}\Delta u = \( \lvert x\rvert^{-\g}\ast \lvert
  u\rvert^2\) u. 
\end{equation}
Note that the nonlinearity $u\mapsto \(
  \lvert x\rvert^{-\g}\ast \lvert u\rvert^2\) u$ is always a smooth
  homogeneous (cubic)
  function of $u$. We assume $2\le \g<\min (4,n)$. A complete scattering
  theory is available in the space $\Sigma$; see
  \cite{GV80,HT87r3}. If we assume moreover $\g>2$, then $\Sigma$ can
  be replaced by $H^1(\R^n)$; see \cite{GVr3energy,Nakanishir3}. 
The counterpart of Lemma~\ref{lem:alg1} is:
\begin{lemma}\label{lem:alg2}
  Let $n\ge 3$ and $2\le \g<\min(4,n)$. Set
  \begin{equation*}
    (q,r)=\( \frac{8}{\g}\virgp \frac{4n}{2n-\g}\).
  \end{equation*}
Then $(q,r)$ is $L^2$-admissible. Set
$\theta =2-4/\g$. 
Then $\theta \in [0,1[$. Define $s=r$ and $k=q/(1-\theta)$. Obviously,
\begin{equation*}
\frac{1}{s}=\frac{1-\theta}{r}+
\frac{\theta}{r}\quad ;\quad  
    \frac{1}{k}=\frac{1-\theta}{q}+
\frac{\theta}{\infty}\virgp
\end{equation*}
and we have $s<\displaystyle
  \frac{2n}{n-\g}$, with
$\displaystyle  \frac{1}{r'}=\frac{1}{r}+\frac{2}{s}
  +\frac{\g}{n}-1$ and $\displaystyle
  \frac{1}{q'}=\frac{1}{q}+\frac{2}{k} $. 
\end{lemma}
We can then proceed as in the energy-subcritical case for the
nonlinear Schr\"odinger equation \eqref{eq:NLS}, in order to prove
Theorem~\ref{theo:r3}. The only difference
is the use of the Hardy--Littlewood--Sobolev inequality. Since the
computations are very similar to those presented in \S\ref{subsec:NLS},
we shall be rather sketchy, and detail only the most important
computation. We set
\begin{equation*}
  F_1=(C\cap L^\infty)(\R;H^1(\R^n))\quad ;\quad F_2 =
  L^q\(\R;W^{1,r}(\R^n)\),
\end{equation*}
where $(q,r)$ is now given by Lemma~\ref{lem:alg2}. It follows from
Strichartz estimates that ($H1$) is satisfied. For $t\in
\overline \R$ and $I$ an interval in $\R$, we have, for $\ell=0$ or $1$:
\begin{align*}
  \Bigg\Vert {\mathbf 1}_{t\in I}&\nabla^\ell \int_{t_0}^t U(t-\tau)\(\(
  \lvert x\rvert^{-\g}\ast (u_1 u_2)\)
  u_3\)(\tau)d\tau\Bigg\rVert_{L^\infty_tL^2_x\cap L^q_tL^r_x} \\
&\lesssim \left\lVert  {\mathbf 1}_{t\in I}\nabla^\ell\(
  \lvert x\rvert^{-\g}\ast (u_1 u_2)\)
  u_3\right\rVert_{L^{q'}_t L^{r'}_x}\\
&\lesssim \Big\lVert  \lVert u_1 \nabla^\ell u_2 \rVert_{L^{s/2}_x}
  \lVert u_3\rVert_{L^r_x}  \Big\rVert_{L^{q'}_t(I)}+ 
\Big\lVert  \lVert u_2 \nabla^\ell u_1 \rVert_{L^{s/2}_x}
  \lVert u_3\rVert_{L^r_x}  \Big\rVert_{L^{q'}_t(I)}\\
&+ \Big\lVert  \lVert u_1 u_2 \rVert_{L^{s/2}_x}
  \lVert \nabla^\ell u_3\rVert_{L^r_x}  \Big\rVert_{L^{q'}_t(I)}\\
&\lesssim \sum_{j=1}^3 \Big\lVert  
  \lVert \nabla^\ell u_j\rVert_{L^r_x}\prod_{j'\not =j}
\lVert u_{j'}\rVert_{L^r_x}  \Big\rVert_{L^{q'}_t(I)} 
\end{align*}
where we have used H\"older and Hardy--Littlewood--Sobolev
inequalities in the space variable. Using H\"older's inequality in
time, we can estimate each term of the above sum by:
\begin{align*}
  \Big\lVert  \nabla^\ell u_j\Big\rVert_{L^q(I;L^r)}&\prod_{j'\not =j}
\lVert u_{j'}\rVert_{L^k(I;L^r)}\\
&\lesssim \Big\lVert  \nabla^\ell
u_j\Big\rVert_{L^q(I;L^r)} \prod_{j'\not =j}\(
\lVert u_{j'}\rVert_{L^q(I;L^r)}^{1-\theta}\lVert
u_{j'}\rVert_{L^\infty(I;L^r)}   \)\\
&\lesssim \prod_{j=1}^3 \lVert {\mathbf 1}_{t\in
  I}u_j\rVert_{F_2}^{1-\theta}\lVert 
u_j\rVert_{F}^{\theta} ,
\end{align*}
where we have used the embedding $H^1\hookrightarrow L^r$. This
estimate suffices to check that Assumption~($H2$) is satisfied (with
$\delta =1-\theta>0$), hence
Theorem~\ref{theo:r3} in the case of $H^1(\R^n)$. In the case of
$\Sigma$ (which allows to consider the value $\gamma=2$), one uses the
operator $J(t)=x+it\nabla$ like in \S\ref{sec:Sigma}, to complete the
proof of Theorem~\ref{theo:r3}.

\subsection{The wave equation}
\label{subsec:NLW}
 We now turn to the case of the nonlinear wave equation
 \begin{equation}
   \label{eq:NLW}
   \d_t^2 u - \Delta u +u^p=0,\quad (t,x)\in \R\times \R^n.
 \end{equation}
In order for the nonlinearity to be analytic, we assume that $p$ is an
integer. Moreover, for the anti-derivative of the nonlinearity to have
a constant sign, we need to assume that $p$ is odd; without this
assumption, scattering for arbitrary large data does not hold. 

The existence of wave and scattering operators in 
 \begin{equation*}
  \Sigma_2 = \{ (f,g)\in H^1(\R^n) \times L^2(\R^n),\quad x\mapsto \lvert
  x\rvert \nabla f(x), x\mapsto \lvert
  x\rvert g(x)\in 
  L^2(\R^n)\}
\end{equation*}
was established in \cite{GV-NLW87}, under the assumption 
\begin{equation*}
  1+\frac{4}{n-1}\le p<1+\frac{4}{n-2}. 
\end{equation*}
As a matter of fact, some values for $p<1+4/(n-1)$ are also allowed
there. See also \cite{BaezSegalZhou} and
\cite{Hidano} for $n=p=3$. With these results, we could certainly
prove that the 
wave and scattering operators are analytic from $\Sigma_2$ to
$\Sigma_2$, for $2\le n\le 4$ and
  \begin{itemize}
\item $p\ge 5$ if $n=2$.
\item $p=3$ or $5$ if $n=3$.
\item $p=3$ if $n=4$. 
  \end{itemize}
We leave out the discussion at this stage, since the estimates based
on the conformal decay are fairly long to write. 
 
The existence of wave and scattering operators in $\dot H^1(\R^n)\times
L^2(\R^n)$ was established in
\cite{BG99,Kapitanski,ShatahStruwe93,ShatahStruwe94} for the
energy-critical case 
\begin{equation*}
  p=1+\frac{4}{n-2},\quad n=3,4.
\end{equation*}
(The space dimensions $3$ and $4$ are the only ones for which the
energy-critical nonlinearity corresponds to an odd integer $p$.) As
stated in Theorem~\ref{theo:NLW}, we shall content ourselves with
these two cases. Note also that from \cite{GV-CMP89}, the existence of
scattering operators in the energy space is known for
energy-subcritical nonlinearities. However, this range for $p$ does
not include odd integers, and we are left with the above two
cases. Also, if we considered only small data scattering, then more
results would be available. We choose not to distinguish too many
cases, and restrict our attention to the framework of
Theorem~\ref{theo:NLW}. 
\bigbreak

Naturally, we have~$D = \dot H^{1}(\R^{n})\times L^{2}(\R^{n})$,
 and
 \begin{equation*}
  F_{1} =( C\cap L^{\infty})(\R; \dot H^{1}(\R^{n}))  
 \times ( C\cap L^{\infty})(\R; L^{2}(\R^{n})).  
 \end{equation*}
As in the case of the
 Schr\"odinger equations 
studied above, the space~$F_{2}$ is defined using Strichartz
 estimates: we set
 \begin{equation*}
   F_2 =\left\{
     \begin{aligned}
       &L^5\(\R;L^{10}(\R^3)\)\times
       L^\infty\(\R;L^2(\R^3)\)&&\quad\text{if }n=3,\\
&L^3\(\R;L^6(\R^4)\)\times L^\infty\(\R;L^2(\R^4)\)&&\quad\text{if
       }n=4. 
     \end{aligned}
\right.
 \end{equation*}
Recall that for $n\ge 3$, and $(q,r)$ satisfying
\begin{equation*}
  \frac{1}{q}+\frac{n}{r}= \frac{n}{2}-1,\quad  6\le r<\infty\text{ if
  }n=3,\ \frac{2n}{n-2}\le r\le
  \frac{2n+2}{n-3}\text{ if
  }n\ge 4,
\end{equation*}
Strichartz estimates yield (see e.g. \cite{GVStriOndes,KT})
\begin{equation*}
  \begin{aligned}
    \left\lVert u\right\rVert_{L^{q}(I;L^{r})} + &
\left\lVert u\right\rVert_{L^{\infty}(I;\dot H^1)} +
\left\lVert \d_t u\right\rVert_{L^{\infty}(I;L^{2})}\\
&\le
 C_{r}\(\left\lVert u_{\mid t=0}\right\rVert_{\dot H^1} +  
 \left\lVert \d_t u_{\mid t=0}\right\rVert_{L^2} + \left\lVert
    \(\d_t^2-\Delta\)u\right\rVert_{L^{1}(I;L^{2})}\),   
  \end{aligned}
\end{equation*}
for some constant $C_r$ independent of the time interval $I$. Note
that the pairs $(5,10)$ and $(3,6)$ are admissible for $n=3$ and
$n=4$, respectively. 

In the case $n=3$, and in view of Example~\ref{ex:ondes}, it is enough
to control $\lVert u_1 u_2 u_3 u_4 u_5\rVert_{L^1(I;L^2)}$ by the
product of the $\lVert u_j\rVert_{L^5(I;L^{10})}$, to verify
Assumption~($H2$). Such an estimate if of course trivially
satisfied. Similarly, for $n=4$, $\lVert u_1 u_2 u_3
\rVert_{L^1(I;L^2)}$ is controlled by the  
product of the $\lVert u_j\rVert_{L^3(I;L^{6})}$. Therefore,
Theorem~\ref{theo:NLW} follows from Lemmas~\ref{lem:meta} and
\ref{lem:scattering}. 

\subsection{The Klein--Gordon equation}
\label{subsec:KG}
We conclude with the case of the Klein--Gordon equation
\begin{equation}
   \label{eq:NLKG}
   \d_t^2 u - \Delta u +u +u^p=0,\quad (t,x)\in \R\times \R^n.
 \end{equation}
As above, we assume that $p$ is an odd integer. The natural energy
space is $D=H^1(\R^n)\times L^2(\R^n)$. For $n\ge 3$, scattering in
the energy space was established in \cite{Brenner84} for
\begin{equation*}
  1+\frac{4}{n}<p\le 1+\frac{4}{n-1},
\end{equation*}
and in \cite{GV85b} for
\begin{equation*}
  1+\frac{4}{n}<p< 1+\frac{4}{n-2}.
\end{equation*}
The case of the low dimensions $n=1$ or $2$ was treated by
K.~Nakanishi \cite{Nakanishi99} (see also \cite{Nakanishi01}), for
$p>1+4/n$. The existence of wave and scattering operators in the
energy-critical case $p=1+4/(n-2)$ in space dimension $n\ge 3$ was
established in \cite{NakanishiIMRN99}. All in all, scattering in the
energy space is known for 
$p>1+4/n$, and $p\le 1+4/(n-2)$ when $n\ge 3$. Such values for $p$
corresponding to an odd integer are exactly those considered in
Theorem~\ref{theo:NLKG}. 
\smallbreak

As pointed out in \cite{GV-NLW87},  this numerology is the
same as in the case of the nonlinear Schr\"odinger equation
\eqref{eq:NLS}. The proof of Theorem~\ref{theo:NLKG} follows
essentially the same lines as the proof of Theorem~\ref{theo:NLS}, up
to the following adaptation. For the space $F_1$, we keep 
\begin{equation*}
  F_1 = \(C\cap L^\infty\) \(\R;H^1(\R^n)\).
\end{equation*}
For the space $F_2$, Sobolev spaces are replaced by Besov spaces:
\begin{align*}
  &F_2 = L^{\g_0} \(\R;B^{1/2}_{\g_0,2}(\R^n)\)\cap L^{\g_1}\(\R\times
  \R^n\),\\
\text{with }& \g_0=2+\frac{4}{n}\text{ and
  }\frac{p-1}{\g_1}+\frac{1}{\g_0} = \frac{1}{\g_0'}. 
\end{align*}
Equation~(3.9) in \cite{Nakanishi99} yields the analogue of the
estimate recalled in Remark~\ref{rem:NLSBesov}:
\begin{align*}
  \lVert u \rVert_{L^\infty (\R;H^1)} +& \lVert \d_t
  u\rVert_{L^\infty(\R;L^2)} + \lVert u \rVert_{F_2}\\
 \le C&\(
  \lVert u_0 \rVert_{H^1} + \lVert u_1 \rVert_{L^2} + \Big\lVert
  \d_t^2 
  u -\Delta 
  u+u \Big\rVert_{L^{\g_0'}\big(\R; B^{1/2}_{\g_0',2}(\R^n)\big)}\),
\end{align*}
The proof of Theorem~\ref{theo:NLKG} then follows the same lines as
the proof of Theorem~\ref{theo:NLS}, up to the technical modifications
which can be found in \cite{Nakanishi99}.

\section{Some consequences}
\label{sec:consequences}

\subsection{Invariant skew-symmetric forms}
\label{sec:form}

Let 
\begin{equation}
  \label{eq:formW}
  \om_{\rm wave}\(u_1,u_2\)(t) := \int_{\R^n}
 \( u_1\d_t u_2-u_2 \d_t u_1\)(t,x)dx.
\end{equation}
It is proved in \cite{MorawetzStrauss} that for the cubic
three-dimensional Klein--Gordon equation (Eq.~\eqref{eq:NLKG} with
$n=p=3$), $ \om_{\rm wave}$ induces a skew-symmetric differential form
on some space $F$ (based on the energy space), which is invariant under
$S$. In \cite{BaezZhou90}, the space $F$ was replaced by the energy
space, in the small data case.  Following the proof of
\cite{MorawetzStrauss}, we have the following extension:
\begin{proposition}
For $m\ge 0$, consider the equation (wave or Klein--Gordon)
\begin{equation*}
  \d_t^2 u-\Delta u +m^2 u +u^p=0.
\end{equation*}
Then under the algebraic assumptions of Theorem~\ref{theo:NLW} (case
$m=0$) or Theorem~\ref{theo:NLKG} (case $m>0$), $\om_{\rm wave}$ induces a
skew-symmetric differential form on the energy space, which is invariant under
$S$. 
\end{proposition}
\begin{proof}[Sketch of the proof]
  Since the proof follows the same lines as in \cite{MorawetzStrauss},
  we shall simply recall the main steps. At least for smooth
  solutions, we compute
  \begin{equation*}
    \frac{d}{dt}\om_{\rm wave} (u_1,u_2) = \int_{\R^n}\(u_2 u_1^p -
    u_1 u_2^p\)dx.
  \end{equation*}
If $u_1$, $u_2$ and $u_3$ solve the above equation, then using the
above relation and expanding
\begin{equation*}
\om_{\rm wave}  \(u_2-u_1,u_3-u_1\)=\om_{\rm wave} \(u_2,u_3\)
+\om_{\rm wave} \(u_1,u_2\) - \om_{\rm wave}  \(u_1,u_3\),  
\end{equation*}
we find
\begin{equation}
  \label{eq:MSW}
  \begin{aligned}
  \frac{d}{dt}\om_{\rm wave}\(u_2-u_1,u_3-u_1\) &= \int\(
  \(u_2^{p-1}-u_3^{p-1}\) (u_2-u_1)u_3\)dx \\
&+ \int\(
  \(u_2^{p-1}-u_1^{p-1}\) (u_3-u_2)u_1\)dx.
\end{aligned}
\end{equation}
Elementary computations show that $(u_1-u_2)(u_1-u_3)(u_2-u_3)$ can be
factored out in the above expression. Now let $u_-,v_-$ and $w_-$ be
in the energy space (whose definition varies whether
$m=0$ or $m>0$). In \eqref{eq:MSW}, we consider $u_1$, $u_2$ and $u_3$
with asymptotic states as $t\to -\infty$ given by $u_-$, $u_-+\eps
v_-$ and $u_-+\eps w_-$, respectively. 
The results of Section~\ref{sec:abstract} show that the image of $v_-$
under $dS(u_-)$ is $v_+$, which is the asymptotic state as $t\to
+\infty$ of $v$, satisfying
\begin{equation*}
  \d_t^2 v-\Delta v +m^2 v +pu^{p-1}v=0,
\end{equation*}
with asymptotic state $v_-$ as $t\to -\infty$ ($v_+=v_-$ if $u\equiv
0$: $S$ is almost the identity near the origin; $v_+$ is implicit
otherwise, see \S\ref{subsec:setting}). 
Integrating \eqref{eq:MSW} over
all $t$, we get:
\begin{equation*}
  \om_{\rm wave} \( (u_2-u_1)_+, (u_3-u_1)_+\) -\om_{\rm wave} \(\eps
  v_-,\eps w_-\) = \O\(\eps^3\), 
\end{equation*}
from the factorization mentioned above. Simplifying by $\eps^2$, the
result follows by letting $\eps \to 0$. 
\end{proof}
In the case of the Schr\"odinger operator, introduce
\begin{align*}
 \om_{\text{Schr\"od}} \(u_1,u_2\)(t)= \IM \int_{\R^n}\(\overline u_1
  u_2\) (t,x)dx.
\end{align*}
Like above, if $u_1$ and $u_2$ solve
\begin{equation*}
  i\d_t u_j+\frac{1}{2}\Delta u_j = F_j,
\end{equation*}
then we have:
\begin{equation*}
  \frac{d}{dt}\om_{\text{Schr\"od}} \(u_1,u_2\)= \RE \int_{\R^n}
  \(\overline F_1 u_2 - 
  \overline u_1 F_2\). 
\end{equation*}
If $u_1$, $u_2$ and $u_3$ solve \eqref{eq:NLS},
we find:
\begin{align*}
  \frac{d}{dt}\om_{\text{Schr\"od}} \(u_2-u_1,u_3-u_1\) &= \int\(
  \lvert u_2\rvert^{p-1} - 
  \lvert u_3\rvert^{p-1} \)\RE (u_2-u_1)\overline u_3\\
& + \int\(\lvert
  u_2\rvert^{p-1} - 
  \lvert u_1\rvert^{p-1} \)\RE (u_3-u_2)\overline u_1 .
\end{align*}
Viewing the right hand side as a polynomial in three unknowns $u_1$,
$u_2$ and $u_3$, we note that it is zero for $u_1=u_2$, $u_3=u_1$ and
$u_2=u_3$. We can then use the same argument as above, to claim that
it yields a contribution of order $\O(\eps^3)$. Proceeding as above,
we have:
\begin{proposition}
  Consider the equation 
\begin{equation*}
  i\d_t u+\frac{1}{2}\Delta u = \lvert u\rvert^{p-1}u. 
\end{equation*}
Under the algebraic assumptions of Theorem~\ref{theo:NLS},
$\om_{\text{Schr\"od}}$ induces a 
skew-symmetric differential form on $H^1(\R^n)$ (or $\Sigma$), which
is invariant under $S$, the scattering operator associated to the
above equation. 
\end{proposition}
Finally, if $u_1$, $u_2$ and $u_3$ solve
\begin{equation*}
  i\d_t u_j+\frac{1}{2}\Delta u_j = \(V\ast \lvert u_j\rvert^2\)u_j,
\end{equation*}
then we find
\begin{align*}
  \frac{d}{dt}\om_{\text{Schr\"od}} \(u_2-u_1,u_3-u_1\) &= \int\(
  V\ast\(\lvert u_2\rvert^{2} - 
  \lvert u_3\rvert^{2}\) \)\RE (u_2-u_1)\overline u_3\\
& + \int\(V\ast\(\lvert
  u_2\rvert^{2} - 
  \lvert u_1\rvert^{2}\) \)\RE (u_3-u_2)\overline u_1 .
\end{align*}
\begin{proposition}
  Consider the equation 
\begin{equation*}
  i\d_t u+\frac{1}{2}\Delta u = \(\lvert x\rvert^{-\g}\ast \lvert
  u\rvert^{2}\)u.  
\end{equation*}
Under the algebraic assumptions of Theorem~\ref{theo:r3},
$\om_{\text{Schr\"od}}$ induces a 
skew-symmetric differential form on $H^1(\R^n)$ (or $\Sigma$), which
is invariant under $S$, the scattering operator associated to the
above equation. 
\end{proposition}

\subsection{Infinitely many conserved quantities}
\label{sec:infiniteconsquant}
In \cite{Baez90,BaezZhou90}, the authors consider the Klein-Gordon equations~(\ref{eq:nklgintro})
with~$p=3$, and prove that the analyticity of the scattering operator (which at the time was only known for
small data) implies the existence of a complete set of   conserved quantities with vanishing Poisson brackets. The proof of~\cite{BaezZhou90} relies upon the construction of invariant skew-symmetric forms, as   in the previous section. Once the form~$ \om_{\rm wave}$ is known, one can construct explicitly a complete
set of integrals of motion~$F_j$, with vanishing Poisson brackets. The statement is given below, in all the cases studied in the paper. We refer to~\cite{BaezZhou90} for the proof of the  result, which can be directly
adapted to the skew-symmetric form~$\om_{\text{Schr\"od}}$.
\begin{proposition}
For each of the equations~(\ref{eq:NLSintro}) to~(\ref{eq:nklgintro}) 
 considered in this paper, and under the algebraic assumptions
of Theorems~\ref{theo:NLS} to~\ref{theo:NLKG}  respectively there is a   family~$F_j$ of analytic functionals acting 
from the space of initial data~$D$ into~$\R$, invariant under the nonlinear evolution, and such that there
is a vector field~$v_j$ in~$D$ such that
$$
dF_j = \om(v_j,\cdot)
$$
where~$\om$ denotes respectively~$\om_{\text{Schr\"od}}$ and~$ \om_{\rm wave}$. Moreover, generically in~$u$,  
for any couple of vector fields~$(v,w)$ in~$T_u D$ such that~$dF_j v = dF_j w = 0$, we  have~$ \om(v,w) = 0
$.
\end{proposition}

This   result can be understood as the existence of a Birkhoff normal form (see e.g. \cite{B03,Grebert} 
for a general definition and a presentation of results). However, for nonlinear equations, Birkhoff normal 
forms are usually employed to establish long time existence results (see 
e.g.  \cite{Bourgain04b,BDGS07}), whereas in our case, they come as a consequence of asymptotic properties of solutions which are already known to exist globally.

\subsection{Inverse scattering}
\label{sec:scattinv}

As noticed in \cite[Theorem~2]{MorawetzStrauss}, knowing the
scattering operator near the origin for a nonlinear equation with
analytic nonlinearity suffices to determine the nonlinearity, since
the coefficients of its Taylor series can be computed by
induction. 
\smallbreak

In \cite{Sasaki}, the first term of the asymptotic expansion of
the scattering operator is shown to fully determine a nonlocal
nonlinearity whose form is known in advance (Hartree type
nonlinearity). This approach is applied in the Schr\"odinger case, as
well as in the Klein--Gordon case. In that case, the nonlinearity need
not be analytic, and only the first nontrivial term of the asymptotic
expansion of $S$ near the origin is needed. Typically, in the same
spirit, consider the nonlinear Schr\"odinger equation
\begin{equation}\label{eq:NLSgen}
  i\d_t u + \frac{1}{2}\Delta u = \lambda \lvert u\rvert^{p-1} u,
\end{equation}
with $\lambda\in \R$ (possibly negative), $p\ge 1+4/n$ and $p\le 1+4/(n-2)$
if $n\ge 3$, not necessarily an integer. For small data, solutions to
\eqref{eq:NLSgen} are global in time, and admit scattering states. To
see this, recall that the nonlinearity in \eqref{eq:NLSgen} is
$H^s$-critical, with
\begin{equation*}
  s = \frac{n}{2} - \frac{2}{p-1}\ge 0. 
\end{equation*}
In the small data  case, Strichartz and Sobolev inequalities show that
global existence and scattering follow from a simple bootstrap
argument (see e.g. \cite{CW89} in the case of $s=0$, \cite{CW90} in
the case $s>0$).  In addition, we have
\begin{equation*}
  W_\pm \(\eps \phi\) =\eps \phi +i\lambda\eps^{p}\int_0^{\pm
  \infty}e^{-i\frac{t}{2}\Delta}\(\left\lvert 
  e^{i\frac{t}{2}\Delta}\phi\right\rvert^{p-1}e^{i\frac{t}{2}\Delta}\phi\)
dt + \O_{H^s}\(\eps^{2p-1}\),
\end{equation*}
hence
\begin{equation*}
  S \(\eps \phi\) =\eps \phi -i\lambda\eps^{p}\int_{-\infty}^{+
  \infty}e^{-i\frac{t}{2}\Delta}\(\left\lvert 
  e^{i\frac{t}{2}\Delta}\phi\right\rvert^{p-1}e^{i\frac{t}{2}\Delta}\phi\)
dt + \O_{H^s}\(\eps^{2p-1}\).
\end{equation*}
See \cite{COMRL} for the proof in the case $s=0$. The proof for $s>0$
follow the same lines, up to the modifications which can be found in
\cite{CW90}. Loosely speaking, the leading order term of
$S(\eps\phi)-\eps \phi$
suffices to determine $\lambda$ and $p$. For instance, 
\begin{equation*}
  p = \lim_{\eps \to 0} \frac{\log \lVert S(\eps \phi) -\eps
  \phi\rVert_{H^s}}{\log \eps},
\end{equation*}
for $\phi$ a Gaussian function, so that the term in $\eps^p$ cannot be
zero.

\subsection{On the complete integrability}
\label{sec:integrable}

When speaking of complete integrability, one has to be rather
cautious: several notions are present in the
literature \cite{Audin2,Integrability}. The weakest definition (which is in fact  useful mainly in a finite dimensional situation)
consists in saying that there 
exists as many conserved quantities as the number of degrees of freedom (infinitely many in infinite dimensional situations), with vanishing Poisson
brackets; this corresponds to the discussion in Section~\ref{sec:infiniteconsquant} above. One can observe that those conserved quantities  may not be relevant in terms of Sobolev norms (see for example~\cite{Bourgain96}).
 In the Hamiltonian case, the  
quantities are the Hamiltonian and first integrals; 
 see e.g. \cite{AM,AudinSMF,AudinBBK}. At a higher (in the  infinite dimensional case)
level of precision, there may exist a nonlinear change of variables
which makes the original equation linear. This is typically the case of
one-dimensional Schr\"odinger equations with cubic nonlinearity
\cite{ZS}, and is related to the existence of Lax pairs \cite{Lax68}.
The strongest notion of integrability consists in trivializing the
equation on some Lie algebra; see e.g. \cite{Helein07}. 

\smallbreak

\subsection*{Acknowledgments}
The authors are grateful to Prof. Tohru Ozawa for pointing out several
references, and to Satoshi Masaki for an early view of his result
\cite{Masaki07}. They also thank Fr\'ed\'eric H\'elein for useful explanations on the notion of complete integrability.
\bibliographystyle{amsplain}
\bibliography{analytic}

\end{document}